\documentclass[preprint]{elsarticle}

\usepackage[hidelinks]{hyperref}

\usepackage[utf8]{inputenc}
\usepackage[T1]{fontenc}

\usepackage[DIV11]{typearea}

\usepackage{amsmath,amsthm}
\usepackage{mathtools}

\usepackage{tikz}
\usepackage{environ} 
\makeatletter
\newsavebox{\measure@tikzpicture}
\NewEnviron{scaletikzpicturetowidth}[1]{%
	\def\tikz@width{#1}%
	\def\tikzscale{1}\begin{lrbox}{\measure@tikzpicture}%
		\BODY
	\end{lrbox}%
	\pgfmathparse{#1/\wd\measure@tikzpicture}%
	\edef\tikzscale{\pgfmathresult}%
	\BODY
}
\makeatother
\usetikzlibrary{shapes}

\usepackage{tabularx}
\usepackage{booktabs}
\usepackage{color}
\usepackage{threeparttable}

\usepackage[markup=nocolor]{changes}
\usepackage{textcomp}
\usepackage{enumitem}
\usepackage{todonotes}

\usepackage{natbib}
\biboptions{round,authoryear}

\newtheorem{thm}{Theorem}

\newtheorem{model}{Model}

\theoremstyle{definition}

\newtheorem{example}[thm]{Example}
\theoremstyle{remark}

\makeatletter
\@ifundefined{c@rownum}{%
	\let\c@rownum\rownum
}{}
\@ifundefined{therownum}{%
	\def\therownum{\@arabic\rownum}%
}{}
\makeatother

 \makeatletter
 \def\user@resume{resume}
 \def\user@intermezzo{intermezzo}
 \newcounter{previousequation}
 \newcounter{lastsubequation}
 \newcounter{savedparentequation}
 \renewenvironment{subequations}[1][]{%
 	\def\user@decides{#1}%
 	\setcounter{previousequation}{\value{equation}}%
 	\ifx\user@decides\user@resume 
 	\setcounter{equation}{\value{savedparentequation}}%
 	\else  
 	\ifx\user@decides\user@intermezzo
 	\refstepcounter{equation}%
 	\else
 	\setcounter{lastsubequation}{0}%
 	\refstepcounter{equation}%
 	\fi\fi
 	\protected@edef\theHparentequation{%
 		\@ifundefined {theHequation}\theequation \theHequation}%
 	\protected@edef\theparentequation{\theequation}%
 	\setcounter{parentequation}{\value{equation}}%
 	\ifx\user@decides\user@resume 
 	\setcounter{equation}{\value{lastsubequation}}%
 	\else
 	\setcounter{equation}{0}%
 	\fi
 	\def\theequation  {\theparentequation  \alph{equation}}%
 	\def\theHequation {\theHparentequation \alph{equation}}%
 	\ignorespaces
 }{%
 	\ifx\user@decides\user@resume
 	\setcounter{lastsubequation}{\value{equation}}%
 	\setcounter{equation}{\value{previousequation}}%
 	\else
 	\ifx\user@decides\user@intermezzo
 	\setcounter{equation}{\value{parentequation}}%
 	\else
 	\setcounter{lastsubequation}{\value{equation}}%
 	\setcounter{savedparentequation}{\value{parentequation}}%
 	\setcounter{equation}{\value{parentequation}}%
 	\fi\fi
 	\ignorespacesafterend
 }
 \makeatother

\newcommand{\rev}[1]{\textcolor{blue}{#1}}

\begin{document}
	\begin{frontmatter}
		\author[1]{Daniela Gaul\fnref{fn2}}
		\fntext[fn2]{Corresponding author.}
		\ead{gaul@math.uni-wuppertal.de}
		\author[1]{Kathrin Klamroth}
		\ead{klamroth@math.uni-wuppertal.de}
		\author[1]{Michael~Stiglmayr}
		\ead{stiglmayr@math.uni-wuppertal.de}
		\title{Event-based MILP models for ridepooling applications}
		\address[1]{School of Mathematics and Natural Sciences, University of Wuppertal, 42119 Wuppertal, Germany}
		
		\journal{European Journal of Operational Research}

\begin{abstract} 
Ridepooling services require efficient optimization algorithms to simultaneously plan routes and pool users in shared rides.
We consider a static dial-a-ride problem (DARP) where a series of origin-destination requests have to be assigned to routes of a fleet of vehicles. Thereby, all requests have associated time windows for pick-up and delivery, and may be denied if they can not be serviced in reasonable time or at reasonable cost. Rather than using a spatial representation of the transportation network we suggest an event-based formulation of the problem, resulting in significantly improved computational times. 
While the corresponding MILP formulations  require more variables than standard models, they have the advantage that capacity, pairing and precedence constraints are handled implicitly.
The approach is tested and validated using a standard IP-solver on benchmark data from the literature. Moreover, the impact of, and the trade-off between, different optimization goals is evaluated on a case study in the city of Wuppertal (Germany).
\end{abstract}

\end{frontmatter}

\section{Introduction}\label{sec:intro}
\noindent

Ridepooling bridges a gap between rather inflexible line-based public transport (e.g.\ buses and trains) and relatively expensive taxi services in the portfolio of transportation services.
A convenient alternative to public transport, but more affordable than a taxi, is the concept of ridepooling: a taxi-like service, typically operated by mini-buses, where users submit pick-up and drop-off locations and a desired pick-up or drop-off time. Users with similar origin or destination are assigned to the same ride whenever this is economically or ecologically useful. Since rides may be shared, users may have to accept longer ride times. In exchange, they can be offered lower fares.
Furthermore, ridepooling services can be easily adapted to varying passenger demands and can thus provide an acceptable service quality also in suburban areas.

To encourage people to use ridepooling services, an efficient and fair route planning is of central importance. From a modelling perspective,  
the \emph{dial-a-ride problem} (DARP) consists of scheduling transportation requests with given pick-up and delivery times  and assigning them to a set of vehicle routes. Common optimization objectives are to accept as many user requests as possible, to minimize the overall routing costs, and to minimize  the maximum and/or the total waiting times and ride times. 
The DARP is motivated by a variety of real-life applications. 
Dial-a-ride systems originated as door-to-door transportation, primarily in rural areas, where public transport is rarely available with buses going only a couple of times a day.
Another application occurs in health care transport: elderly, injured or disabled persons are often not able to drive on their own or to use public transport, and dial-a-ride services are a convenient alternative to expensive taxi services. Recently, several so-called `on-demand ridepooling services' have emerged. Prominent examples are Uber\footnote{\url{https://www.uber.com/de/en/ride/uberpool/}},
DiDi Chuxing\footnote{\url{https://www.didiglobal.com/travel-service/taxi}}, Lyft\footnote{\url{https://www.lyft.com/rider}} or moia\footnote{\url{https://www.moia.io/en/how-it-works}}. Ridepooling services complement public transport  and are usually established in urban areas. They are an alternative to motorized private transport and thus help to reduce the number of vehicles in the cities.
This paper is motivated by a particular ridepooling service that was established in the city of  Wuppertal (Germany) in collaboration with the project bergisch.smart\footnote{\url{https://www.bergischsmartmobility.de/en/the-project/}}. This ridepooling service is run by the local public transport provider and is designed to complement the bus service.

In this paper, we focus on the \emph{static} DARP assuming that all user requests are known in advance, i.e.\ 
before the DARP is solved. This is in contrast to the \emph{dynamic} DARP where requests are revealed during the day and vehicle routes are updated whenever a new request arrives. 
While the static DARP is important on its own right and most research focuses on the static DARP, see e.\,g.\ \cite{Ho2018}, we note that static DARP models can also be extended to the dynamic scenario by using a rolling-horizon strategy, see e.g.\ \cite{Gaul2021A}.

\subsection{Literature Review}\label{sec:lit}
\noindent
The following literature review focuses specifically on variants of the static multi-vehicle dial-a-ride problem with time windows. For a broader review the reader is referred to \cite{Cordeau2007}  who cover research on the DARP up to 2007, and to \cite{Molenbruch2017} or \cite{Ho2018} for more recent surveys. 

Early work  is carried out by \cite{Jaw1986}, who develop one of the first heuristics for the multi-vehicle DARP. Depending on the users' earliest possible pick-up times, the heuristic determines the cheapest insertion position in an existing route in terms of user satisfaction and operator costs.  
\cite{Desrosiers1991}, \cite{Dumas1989} and \cite{Ioachim1995} first identify groups of users to be served within the same area and time. In a second step, these `clusters' are combined to obtain feasible vehicle routes.
This decomposition approach, however, leads in general to suboptimal solutions.
\cite{Cordeau2003} are the first to apply tabu search to the DARP. In doing so they make use of a simple neighborhood generation technique by moving requests from one route to another. Their reported computational results on real-life data as well as on randomly generated instances validate the efficiency of the heuristic. Recent work on the static multi-vehicle DARP related to tabu search is based on \cite{Cordeau2003}, see, for example, \cite{Detti2017, GUERRIERO2013, Kirchler2013, Paquette2013}.
\cite{Cordeau2006} is, to the best of our knowledge, the first to apply a branch-and-cut (B\&C) algorithm to the DARP. In addition to valid inequalities derived from the traveling salesman problem, the vehicle routing problem, and the pick-up and delivery problem, new valid inequalities are added as cuts to a mixed-integer linear programming (MILP) formulation of the DARP. While the MILP model of \cite{Cordeau2006} is based on three-index variables, and is a standard formulation of the basic DARP that has been used as well for further problem extensions (see \cite{Ho2018}), \cite{Ropke2007} propose an alternative MILP formulation of the DARP using two-index variables and derive three more classes of valid inequalities and combine them with some of the cuts presented by \cite{Cordeau2006} to develop a B\&C algorithm. Some of these inequalities  are employed in more recent works on more complex versions of the DARP, see, for example, \cite{Braekers2014, Parragh2011}.  The former work proposes a B\&C algorithm to exactly solve small-sized instances as well as a deteministic annealing meta-heuristic for larger problem instances, while \cite{Parragh2011} proposes  deterministic annealing and insertion heuristics, respectively. 
There exist several exact solution approaches that are based on branch-and-bound (B\&B) frameworks. For example, \cite{Cortes2010} combine B\&C and Benders decomposition to solve a DARP that allows users to swap vehicles at specific locations during a trip. \cite{Liu2015} devise a B\&C algorithm and formulate two different models for a DARP with multiple trips, heterogeneous vehicles and configurable vehicle capacity. They introduce eight families of valid inequalities to strengthen the models. \cite{Qu2015} propose an exact method based on BCP with heterogeneous vehicles and configurable vehicle types. \cite{Gschwind2015} develop a new exact BCP algorithm for the DARP which outperforms the B\&C algorithm suggested by \cite{Ropke2007}. They show that their algorithm is capable of solving an instance with $8$ vehicles and $96$ users which is, to the best of our knowledge, to date the largest instance of the standard DARP solved by an exact method. 

Since exact branching based procedures are usually limited to small to medium-sized instances, there has been some work on hybridizing B\&B-like methods with heuristics. For example, \cite{Hu2014} apply a branch-and-price (B\&P) algorithm to a DARP  with time-dependent travel times, in which the pricing subproblem is solved by large neighborhood search. As a side result, they observe that the length of time windows can have a significant impact on the size of the fleet, the objective function value, the CPU time, the average ride time or the average pick-up time delay. 
While exact solution approaches are often based on B\&B, at the level of heuristics and meta-heuristics there is a wide set of methods that can be applied.  
\cite{Jorgensen2007} propose an algorithm which assigns passengers to vehicles using a genetic algorithm. In a second step, routes are constructed sequentially by means of a nearest neighbor procedure. Genetic algorithms are also used by \cite{Cubillos2009}, who  obtain slightly better results than \cite{Jorgensen2007}. \cite{Atahran2014} devise a genetic algorithm within a multi-objective framework. \cite{Belhaiza2017} and \cite{Belhaiza2019}  combine genetic crossover operators with variable neighborhood search  and adaptive large neighborhood search, respectively. A hybrid genetic algorithm is also developed by \cite{Masmoudi2017}, who according to \cite{Ho2018}, achieve the second-best results on benchmark instances of the standard DARP. 
\cite{Hosni2014} present a MILP formulation which is then decomposed into subproblems using Lagrangian relaxation. Feasible solutions are found by  modifying the subproblems in a way that generates a minimal incremental cost in terms of the objective function. 
\cite{Gschwind2019} develop a meta-heuristic solution method for the DARP, which is at the moment probably the most efficient heuristic method to solve the standard DARP.  They adopt an adaptive large neighborhood method that was developed by \cite{Ropke2006} for the pick-up and delivery problem with time windows. The resulting algorithm  integrates a dynamic programming algorithm into large neighborhood search. \cite{Souza2020} present a simple heuristic based on variable neighborhood search combined with a set covering strategy. A heuristic repair method based on a greedy strategy, that reduces the infeasibility of infeasible intital routes is developed by \cite{Chen2020}. Their repair method improves about 50\% of the fixing ability of a local search operator.

While most of the literature on the DARP focuses on the single-objective DARP, the DARP is examined from a multi-objective optimization perspective as well. As stated by \cite{Ho2018}, this line of research can be divided into three main categories: use of a weighted sum objective, a lexicographic approach and determination of (an approximation of) the Pareto front, the set of all non-dominated solutions. Weighted sum objectives have been used e.\,g.\ by \cite{Jorgensen2007, Melachrinoudis2007, Mauri2009, Kirchler2013} or \cite{Bongiovanni2019}. 
In addition to the customers' total transportation time, which is a common objective in single-objective DARPs, the weighted sum objective introduced by \cite{Jorgensen2007} is composed of the total excess ride time, the customers' waiting time, the drivers' work time as well as several penalty functions for the violation of constraints. The weights on these criteria are chosen with respect to the relative importance of the criteria from a user-perspective. 
While the methods to solve the DARP with a weighted sum objective range from a genetic algorithm (\cite{Jorgensen2007}), to tabu search (\cite{Melachrinoudis2007, Kirchler2013}), simulated annealing (\cite{Mauri2009}) or B\&C (\cite{Bongiovanni2019}), most of the authors either adapt the weights introduced by \cite{Jorgensen2007} or choose them with respect to individual assessments of their relative importance. \cite{Garaix2010, Schilde2014} and \cite{Luo2019} consider a lexicographic ordering of the objective functions w.r.t.\  their relative importance. In this context, \cite{Luo2019} have proposed a two-phase BCP-algorithm and a strong trip-based model. The model is based on a set of nondominated trips that are enumerated by a label-setting algorithm in the first phase, while in the second phase the model is solved using Benders decomposition and BCP. The Pareto front is examined e.\,g.\ by \cite{Parragh2009, Zidi2012, Chevrier2012, Paquette2013, Atahran2014,  Hu2019} or \cite{Viana2019}. A comparison of six evolutionary algorithms to solve the multi-objevtive DARP is provided by \cite{Guerreiro2020}.

\paragraph{Contribution}
In this paper we suggest an event-based formulation of ridepooling problems which is based on an abstract graph model and, in contrast to other approaches, has the advantage of implicitly encoding vehicle capacities as well as pairing and precedence constraints. The event-based graph model is the basis for two alternative mixed-integer linear programming formulations that can be distinguished by the approach towards handling ride time and time window constraints. The MILP formulations are pseudo-compact, that is, the size of the formulations is polynomial w.r.t.\ the number of requests but depends on the number of vehicles. The presented models can be used to solve small- to medium-sized benchmark instances in a short amount of time. In particular, we show that both suggested MILP formulations outperform the standard compact 3-index formulation of \cite{Cordeau2006} and the compact 2-index formulation for the closely related Pick-up and Delivery Problem with Time Windows (PDPTW) of \cite{Furtado2017}, which can be extended to a MILP formulation for DARP by introducing additional constraints.
Furthermore, conflicting optimization goals reflecting, e.g., economic efficiency and customer satisfaction are combined into weighted sum objectives. In addition, we distinguish between average case and worst case formulations. The impact of these objective functions is tested using artificial request data for the city of Wuppertal, Germany.

This paper is organized as follows. After a formal introduction of the problem in Section~\ref{sec:prob}, the event-based graph and based on that two variants of an MILP formulation of the static DARP are presentend in Section~\ref{sec:model}. The models are compared and tested using CPLEX, in Section~\ref{sec:experiments}. The paper is concluded with a summary and some ideas for future research in Section~\ref{sec:concl}.


\section{Problem Description}\label{sec:prob}
\noindent
The DARP considered in this paper is defined as follows: Let $n$ be the number of users submitting a transport request. Each request $i\in R\coloneqq \{1, \ldots, n\}$ specifies a pick-up location $i^+$ and a drop-off location $i^-$. The sets of pick-up and drop-off locations are denoted by $P\coloneqq \{1^+, \ldots, n^+\}$ and  $D\coloneqq \{1^-, \ldots, n^-\}$, respectively. A homogeneous fleet of vehicles  $K$ with capacity $Q$ is situated at the vehicle depot, which is denoted by $0$. A number of requested seats $q_i \geq 1$  and a service duration of $s_i\geq 0$ are associated with each request $i\in R$ and we set $q_{i^+} = q_{i^-}\coloneqq q_i$, $s_{i^+} = s_{i^-}\coloneqq s_i$ and $q_0 \coloneqq 0$ as well as $s_0 \coloneqq 0$.  The earliest and latest time at which service may start at a request location $j\in P\cup D\cup \{0\}$ is given by $e_j$ and $\ell_j$, respectively, so that each request location has an associated time window $\left[ e_j, \ell_j \right]$.  At the same time, the time window $\left[e_0, \ell_0 \right]$ corresponding to the depot contains information about the overall start of service $e_0$ and the maximum duration of service $T\coloneqq \ell_0 - e_0$. The maximum ride time corresponding to request $i \in R $ is denoted by $L_i$.
A feasible solution to the DARP consists of at most $\vert K\vert$ vehicle routes which start and end at the depot. If a user is served by vehicle $k\in K$, the user's pick-up and drop-off location both have to be contained in this order in vehicle $k$'s route. The vehicle capacity of $Q$ may not be exceeded at any time. The start of service at every location has to be within the time windows. It is possible to reach a location earlier than the start of service and wait. At each location $j$, a service duration of $s_j$ minutes is needed for users to enter or leave the vehicle. In addition, the acceptable ride time of each user $i$ is bounded from above by $L_i$. Vehicles have to return to the depot at least $T$ minutes after the time of the overall start of service $e_0$. We note that this is a common setting for the static DARP, see, e.g., \cite{Gschwind2019, Ropke2007}.
Moreover, dial-a-ride problems are often characterized by conflicting objectives. In this paper, we focus on
\begin{itemize}[itemsep=2pt,parsep=2pt]
	\item economic efficiency, e.g., minimization of routing costs, and
	\item user experience, e.g., minimization of unfulfilled user requests, waiting time and ride time, 
\end{itemize} 
and combine these objectives into an overall \emph{weighted sum objective}. Note that this weighted sum objective can be interpreted as a scalarization of a multi-objective model in which all objective functions are equitably considered. Consequently, by variation of the weights every supported efficient solution, that is every efficient solution on the convex hull of the Pareto front, can be determined as a solution of the weighted sum scalarization, see e.g.\ \cite{ehrgott05multicriteria}.
Other options would be to consider a lexicographic multi-objective model or to treat one objective function, e.g., economic efficiency, as objective and to impose constraints on the other to ensure a satisfactory user experience, which corresponds to the \(\varepsilon\)-constraint scalarization.

Mixed integer linear programming formulations based on these assumptions are presented in the following section.

\section{Mixed Integer Linear Programming Formulations}\label{sec:model}
\noindent
Classical models for the DARP are based on a directed graph that can be constructed in a straight-forward way by identifying all pick-up and drop-off locations and the depot by nodes in a complete directed graph, see e.g.\ \cite{Cordeau2006} or \cite{Ropke2007}.
In contrast to this, we suggest an \emph{event-based graph} $G=(V, A)$ in which the node set $V$ consists of  $Q$-tuples representing feasible user allocations together with information on the most recent pickup or drop-off. As will be explained later, this modelling approach has the advantage that many of the complicating constraints can be encoded directly in the network structure: vehicle capacity, pairing (i.e.\ users' pick-up and drop-off locations need to be served by the same vehicle) and precedence constraints (i.e.\ users' pick-up locations need to be reached before their drop-off locations are reached) are implicitly incorporated in the event-based graph. Moreover, a directed arc $a\in A$ is only introduced between a pair of nodes if the corresponding sequence of events is feasible w.r.t.\ the corresponding user allocation. Thus, a tour in $G$ is feasible if it does not violate constraints regarding time windows or ride times. 

Using the event-based graph, the DARP can be modeled as a variant of a minimum cost flow problem with unit flows and with additional constraints. In particular, we consider circulation flows in $G$ and identify the depot with the source and the sink of a minimum cost flow problem. Each dicycle flow in $G$ then represents one vehicle's tour.

\subsection{Event-based graph model}\label{subsec:event}
\noindent
The event-based MILP formulations presented in this paper are motivated by the work of \cite{bertsimas19online} who propose an optimization framework for taxi routing, where only one passenger is transported at a time. Their algorithm can handle more than $25,000$ users per hour. \cite{bertsimas19online} propose a graph-based formulation in which an arc $(i, j)$ represents the decision to serve passenger $j$ directly after dropping off passenger $i$.

The allocation of users to a vehicle with capacity $Q$ can be written as a $Q$-tuple. Unassigned user slots (i.e., empty seats) are indicated by zeros. If, for example, $Q=3$ and two requests, request $1$ and request $2$ with $q_1=q_2=1$, are assigned to the same vehicle, the user allocation may be, for instance, represented by the tuple $\left(2,1,0\right)$. Accordingly, an empty vehicle is represented by $(0,0,0)$. To additionally incorporate information about the most recent pick-up or drop-off location visited, we write
\begin{equation*}
	(2^+,1,0) \quad \text{or} \quad (2^-,1,0),
\end{equation*}
if user $1$ is seated and user $2$ has  just been picked-up or dropped-off, respectively. Note that this encoding implicitly specifies the location of the vehicle since, in this example, user $2$ is picked-up at his or her pick-up location $2^+$, i.e., the $3$-tuple $(2^+,1,0)$ can be associated with the location $2^+$, and user $2$ is dropped-off at his or her drop-off location $2^-$, i.e., the $3$-tuple $(2^-,1,0)$ can be associated with the location $2^-$. 

Since all permutations of such a $Q$-tuple specify the same user allocation, we order  the components of the $Q$-tuple such that the first component contains the information regarding the last pick-up or drop-off stop and the remaining $Q-1$ components are sorted in descending order. Users can only be placed together in a vehicle if the total number of  requested seats does not exceed the vehicle capacity. This constraint limits the possible combinations of users in the vehicle and hence the set of possible $Q$-tuples representing user allocations.  As stated by \cite{Cordeau2006}, time window and ride time constraints can be used to identify incompatible user pairs. For example, given users $i, j \in R$, the following reductions are possible:
\begin{itemize}
	\item if it is infeasible to visit locations $ i^+, i^-, j^+, j^-$ both in the order $ j^+ \rightarrow i^+\rightarrow j^- \rightarrow i^-$ and $ j^+ \rightarrow i^+\rightarrow i^- \rightarrow j^-$,  the set $\{(v_1, v_2, \ldots, v_Q)\in V\colon v_1 = i^+, j\in \{v_2, \ldots, v_Q\}\}$ can be removed from the node set
	\item if it is infeasible to visit locations $ i^+, i^-, j^+, j^-$ both in the order $j^+ \rightarrow i^+ \rightarrow i^- \rightarrow j^-$ and $i^+ \rightarrow j^+ \rightarrow i^- \rightarrow j^-$, the set $\{(v_1, v_2, \ldots, v_Q)\in V\colon v_1 = i^-, j\in \{v_2, \ldots, v_Q\}\}$ can be removed from the node set 
\end{itemize}
The number of possible $Q$-tuples can be reduced significantly by adapting this strategy to identify incompatible user allocations. For this purpose, given requests $i,j\in R$, let $f^1_{i,j}$ and $f^2_{i,j}$ encode the feasibility of the paths $j^+ \rightarrow i^+ \rightarrow j^- \rightarrow i^-$ and $j^+ \rightarrow i^+ \rightarrow i^- \rightarrow j^-$, respectively, w.r.t.\ ride time and time window constraints. To simplify the notation, let $f^1_{i,0} = f^1_{0,i} = f^2_{i,0} = f^2_{0,i} = 1$.

Now the DARP can be represented by a directed graph $G = (V,A)$, where the node set $V$ represents events rather than geographical locations. The set of event nodes corresponds to the set of all feasible user allocations. The set of all nodes that represent an event in which a request (or user) $i\in R$ is picked up is called the set of \emph{pick-up nodes} and is given by 
\begin{equation*}
	\begin{split}
		V_{i^+} \coloneqq \Biggl\lbrace  (v_1, v_2, \ldots, v_Q) \colon  v_1 = i^+,\; &v_j \in R \cup \{0\} \setminus \{i\}, \; f^1_{i,v_j} + f^2_{i,v_j} \geq 1\; \forall j \in\{2,\ldots,Q\},\\
		&\Bigl(v_j > v_{j+1} \vee v_{j+1}=0 \Bigr) \; \forall j \in\{2,\ldots,Q-1\}, \; \sum_{j=1}^Q q_{v_j} \leq Q\Biggr\rbrace.
	\end{split}
\end{equation*}
Similarly, the set of \emph{drop-off nodes} corresponds to events where a request (or user) $i \in R$ is dropped off and is given by
\begin{equation*}
	\begin{split}
		V_{i^-} \coloneqq \Biggl\lbrace  (v_1, v_2, \ldots, v_Q) \colon  v_1 = i^-,\; &v_j \in R \cup \{0\} \setminus \{i\}, \; f^1_{v_j,i} + f^2_{i,v_j} \geq 1\; \forall j \in\{2,\ldots,Q\},\\
		&\Bigl(v_j > v_{j+1} \vee v_{j+1}=0 \Bigr) \; \forall j \in\{2,\ldots,Q-1\}, \; \sum_{j=1}^Q q_{v_j} \leq Q\Biggr\rbrace.
	\end{split}
\end{equation*}
Note that for each request $i \in R$ there is only one pick-up and one drop-off location, but in general more than one potential pick-up node and drop-off node. 
Hence, there is a unique mapping of nodes to locations, while a location may be associated with many different nodes. The depot-node $\boldsymbol{0}\coloneqq (0,\ldots,0)$ is associated with the location of the depot, and we write $V_0 \coloneqq \{\boldsymbol{0}\}$. The overall set of nodes $V$ is then given by 
\begin{equation*}
	V = V_0 \cup \bigcup_{i = 1}^n V_{i^+} \cup \bigcup_{i = 1}^n V_{i^-}. 
\end{equation*}

The arc set $A$ of $G$ is defined by the set of possible transits between pairs of event nodes in $V$. It is composed of six subsets, i.e.,
\[
A = \bigcup_{i = 1}^6 A_i,  
\]
where the subsets $A_i$, $i=1,\dots,6$ are defined as follows:
\begin{itemize}
	\item Arcs that describe the transit from a pick-up node in a set $V_{i^+}$, i.e., from a user $i$'s pick-up location, to a drop-off node in $V_{j^-}$, i.e., to the drop-off location of a user $j$, where $j=i$ is possible, but where $j$ may also be another user from the current passengers in the vehicle:
	\begin{equation*}
		\begin{split}
			A_1 \coloneqq \Bigl\lbrace  \bigl(\bigl(i^+, v_2, \ldots, v_Q\bigr), \bigl(j^-, w_2, \ldots, w_Q\bigr)\bigr)\in V\times V   \colon   \{j, w_2, \ldots, w_Q\}=\{i, v_2, \ldots, v_Q\}\Bigr\rbrace.
		\end{split}
	\end{equation*} 
	\item Arcs that describe the transit from a pick-up node in a set $V_{i^+}$, i.e., from a user  $i$'s pick-up location, to another pick-up node from a set $V_{j^+}$ with $j\neq i$, i.e., to another user $j$'s pick-up location:
	\begin{equation*}
		\begin{split}
			A_2 \coloneqq \Bigl\lbrace  \bigl(\bigl( i^+, v_2, \ldots, v_{Q-1}, 0\bigr) , \bigl(j^+, w_2, \ldots, w_Q\bigr)\bigr)\in V\times V  \colon \{i, v_2, \ldots, v_{Q-1}\}=\{w_2, \ldots, w_Q\} \Bigr\rbrace.
		\end{split}
	\end{equation*}
	\item Arcs that describe the transit from a drop-off node in a set $V_{i^-}$, i.e., from a user $i$'s drop-off location, to a pick-up node in a set $V_{j^+}$, $j\neq i$, i.e., to another user  $j$'s pick-up location:
	\begin{equation*}
		A_3 \coloneqq \Bigl\lbrace \bigl(\bigl(i^-, v_2, \ldots, v_Q\bigr) , \bigl( j^+, v_2, \ldots, v_Q\bigr)\bigr) \in V\times V \colon  i \neq j 
		\Bigr \rbrace.
	\end{equation*}
	\item Arcs that describe the transit from a drop-off node in a set $V_{i^-}$, i.e., from a user $i$'s drop-off location, to a node in $V_{j^-}$, $j\neq i$, i.e., to another user $j$'s drop-off location:
	\begin{equation*}
		\begin{split}
			A_4 \coloneqq \Bigl\lbrace \bigl(\bigl( i^-, v_2, \ldots, v_Q\bigr) , \bigl( j^-, w_2, \ldots, w_{Q-1}, 0\bigr)\bigr)\in V\times V \colon  \{ v_2, \ldots, v_Q\}=\{ j, w_2, \ldots, w_{Q-1}\} 
			\Bigr\rbrace.
		\end{split}
	\end{equation*}
	\item Arcs that describe the transit from a drop-off node in a set $V_{i^-}$, i.e., from a  user $i$'s drop-off location, to the depot:
	\begin{equation*}
		A_5 \coloneqq \Bigl\lbrace \bigl(\bigl( i^-, 0, \ldots,0\bigr) , (0,\ldots,0)\bigr)\in V\times V \Bigl\rbrace .
	\end{equation*}
	\item Arcs that describe the transit from the depot to a pick-up node in a set $V_{i^+}$, i.e., to a user $i$'s pick-up location:
	\begin{equation*}
		A_6 \coloneqq \Bigl\lbrace \bigl((0,\ldots,0) , \bigl(i^+,0,\ldots,0\bigr)\bigr)\in V\times V \Bigr\rbrace. 
	\end{equation*}
\end{itemize}

\begin{example}
	We give an example of the event-based graph $G=(V,A)$ with three users and vehicle capacity $Q=3$. Let $R=\{1,2,3\}$, $q_1 = q_2 = 1$ and $q_3 = 3$. We omit time windows, travel and ride times in this example. By the above definitions we obtain the graph illustrated in Figure~\ref{fig:example}. Note that there are no nodes $v\in V$ that simultaneously contain users $1$ (i.e., $1^+$ or $1^-$) and $3$ (i.e., $3^+$ or $3^-$) as the total number of requested seats specified by these users exceeds the vehicle capacity. Similary, the seats requested by users $2$ and $3$ together exceed the vehicle capacity of three.  Two feasible tours for a vehicle in $G$ are given, for example, by the dicycles 
	\begin{equation*}
		\begin{split}
			C_1 =\Bigl\lbrace & \big(\big(0,0,0\big),\big(1^+,0,0\big)\big),\ \big(\big(1^+,0,0\big),\big(2^+,1,0\big)\big),\ \big(\big(2^+,1,0\big),\big(2^-,1,0\big)\big), \\  
			& \big(\big(2^-,1,0\big),\big(1^-,0,0\big)\big),\  \big(\big(1^-,0,0\big),\big(0,0,0\big)\big)\Bigr\rbrace
		\end{split}
	\end{equation*}
	and
	\begin{equation*}
		C_2 = \Bigl\lbrace \big(\big(0,0,0\big),\big(3^+,0,0\big)\big),\ \big(\big(3^+,0,0\big), \big(3^-,0,0\big)\big),\ \big(\big(3^-,0,0\big),\big(0,0,0\big)\big)\Bigr\rbrace.
	\end{equation*}

	\begin{figure}
		\begin{scaletikzpicturetowidth}{\textwidth}
			\footnotesize
			\begin{tikzpicture}[scale = \tikzscale]
				\useasboundingbox (-12,-11) rectangle (12,0.75);
				\node (null) at (0,0) {$(0,0,0)$};
				\node  (100) at (-10,-3) {$\left(1^+,0,0\right)$};
				\node (200) at (-2,-3) {$\left(2^+,0,0\right)$};
				\node (300) at (6,-3) {$\left(3^+,0,0\right)$};
				\path (null) edge[-latex] (100);
				\path (null) edge[-latex] (200);
				\path (null) edge[-latex] (300);
				\node (400) at  (-6,-3) {$\left(1^-,0,0\right)$};
				\node (500) at (2,-3) {$\left(2^-,0,0\right)$};
				\node (600) at (10,-3) {$\left(3^-,0,0\right)$};
				
				\path
				(100) edge[-latex] (400)
				(200) edge[-latex] (500)
				(300) edge[-latex] (600);
				
				\path
				(400) edge[-latex] (null)
				(500) edge[-latex] (null)
				(600) edge[-latex] (null);
				
				\path
				(400) edge[-latex] (200)
				(400) edge[-latex, bend angle = 20, bend right] (300);
				
				\path
				(500) edge[-latex, bend angle = 20, bend left] (100)
				(500) edge[-latex] (300);
				
				\path
				(600) edge[-latex, bend angle = 20, bend left] (100)
				(600) edge[-latex, bend angle = 20, bend left] (200);	
				
				\node (210) at (-2,-6) {$\left(2^+,1,0\right)$};
				
				\path
				(100) edge[-latex] (210);
				
				\node (120) at (-10,-6) {$\left(1^+,2,0\right)$};
				
				\path
				(200) edge[-latex] (120);
				
				\node (420) at  (-6,-9) {$\left(1^-,2,0\right)$};
				\node (510) at  (2,-9) {$\left(2^-,1,0\right)$};
				
				\path 
				(210) edge[-latex] (420)
				(210) edge[-latex] (510);
				
				\path 
				(120) edge[-latex] (420)
				(120) edge[-latex] (510);

				\draw (510) edge[-latex, bend angle = 15, bend right] (400);
				\draw (420) edge[-latex, bend angle = 30, bend right] (500);
				
			\end{tikzpicture}
		\end{scaletikzpicturetowidth}
		\caption{Graph representation of an example with three users\label{fig:example}.}
	\end{figure}

\end{example}

In order to evaluate the worst case complexity of this event-based graph representation of the DARP, we first evaluate the respective cardinalities of the node set $V$ and of the arc set $A$. Note that the number of nodes and arcs in the event-based graph model depends on the vehicle capacity, the number of users and the number of requested seats per user and can potentially be reduced based on e.g.\ precedence or time window constraints. Given $Q$ and $n$, the number of nodes and arcs is maximal if all users request only one seat, i.e.\ if $q_i=1$ for all $i\in R$. Thus, we obtain the worst case bound
\begin{equation*}
	\vert V \vert \leq 1  + 2\,n \sum_{j=0}^{Q-1} \binom{n-1}{j}, 
\end{equation*}
where 
\begin{itemize}
	\item ``1'' corresponds to the depot node, 
	\item the factor $2n$ corresponds to the first component of the nodes being equal to $i^+$ or $i^-$ with $i\in R = \{1,\ldots,n\}$,
	\item and the sum corresponds to the allocation of users to the remaining seats $v_2, \ldots, v_Q$, of which $j \in \{0, \ldots, Q-1\}$ are occupied.
\end{itemize} 
For the arc set, we obtain the worst case bound
\begin{equation*}
	\begin{split}
		\vert A\vert \leq \  2\,n + n \sum_{j=0}^{Q-1} \binom{n-1}{j}(j+1)  + 3\,n\,(n-1)\sum_{j=0}^{Q-2} \binom{n-2}{j} + \frac{n (n-1)\cdot\ldots\cdot (n-Q)}{(Q-1)!}.
	\end{split}
\end{equation*} 

This can be seen by counting the number of outgoing arcs for each node with $Q-j+1$ zero entries and adding them up for $j\in \{0, \ldots, Q-1\}$, together with $2n$ arcs corresponding to the depot node.  We use the convention that $\binom{m}{k}\coloneqq 0 $ when $k > m$. From the above formulas, we deduce that the number of nodes is bounded by $\mathcal{O}(n^Q)$ for $n\geq Q$ and the number of arcs is bounded by $\mathcal{O}(n^{Q+1})$ for $n\geq Q+1$. This is in general considerably more than what is obtained when a classical, geometrical DARP model on a complete directed graph is used, which has only $\mathcal{O}(n)$ nodes and $\mathcal{O}(n^2)$ arcs, see e.g.\ \cite{Cordeau2006, Ropke2007}. However, in practice, ridepooling services are usually operated by taxis or mini-busses, so that $Q \in \{3,6\}$. Moreover, the number of nodes and thereby the number of arcs, reduce substantially if we do not consider the ``worst-case-scenario'' $q_i = 1$ for all $i\in R$, in which all combinations of requests are possible user allocations in the vehicle, and if we consider time windows, travel and ride times to determine infeasible user pairs. Besides that, the event-based formulation has the clear advantage that important constraints like vehicle capacity constraints, pairing constraints and precedence constraints that have to be formulated in classical models are implicitly handled using the event-based graph model, as will be seen in the next section.

\subsection{Event-based MILP models}
\noindent 
With the above definitions, the DARP can be modeled as a minimum cost integer flow problem with additional constraints, where both the source and the sink are represented by the node $\boldsymbol{0}$. We will formulate corresponding MILP models in the subsequent subsections. Therefore, we need the following additional parameters and variables, which are also summarized in Tables~\ref{tab:paras} and \ref{tab:vars}.

Since each node in $V$ can be associated with a unique request location $j \in P\cup D \cup \{0\}$, we can associate a routing cost $c_a$ and a travel time $t_a$ with each arc $a=(v,w)\in A$. More precisely, both values $c_a$ and $t_a$ are calculated by evaluating the actual routing cost and travel time from the location associated with $v$ to the location associated with $w$. We assume that all routing costs and all travel times are nonnegative and satisfy the triangle inequality. Finally, let $\delta^{\text{in}}(v)\coloneqq \{(w,v) \in A \}$ and $\delta^{\text{out}}(v)\coloneqq\{(v,w) \in A \}$ denote the set of incoming arcs of $v$ and the set of outgoing arcs of $v$, respectively. 
Moreover, for $a\in A$ let the variable value $x_a = 1$ indicate that arc $a$ is used by a vehicle, and let $x_a = 0$ otherwise. Thus, a vehicle tour is represented by a sequence of events in a dicycle $C$ in $G$ where $x_a=1$ for all $a \in C$. A request is matched with a vehicle if the vehicle's tour, i.e., the sequence of events in the corresponding dicycle, contains the event of picking-up and dropping-off of the corresponding user. Since we allow users' requests to be denied, let the variable value $p_i = 1$ indicate that request $i\in R$ is accepted, and let $p_i = 0$ otherwise. Let the variable $B_v$ store the information on the beginning of service time at node $v \in V$. Recall that the beginning of service time has to be within the associated time window of the respective location of node $v$.

\begin{table}
	\centering
	\begin{threeparttable}
		\small
		\begin{tabularx}{0.8\textwidth}{cX}
			\toprule
			Parameter & Description \\
			\midrule
			$n$ & number of  transport requests \\
			$R$ & set of transport requests \\
			$i^+$, $i^-$ & pick-up and drop-off location of request $i$ \\
			$P$, $D$ & set of pick-up locations/requests, set of drop-off locations\\
			$K$ & fleet of vehicles \\
			$Q$ & vehicle capacity \\
			$q_j$ & load associated with location $j$ \\
			$s_j$ & service duration associated with location $j$ \\
			$\left[e_j, \ell_j \right]$ & time window associated with location $j$ \\
			$T$ & maximum duration of service \\
			$L_i$ & maximum ride time associated with request $i$ \\
			$V$ & node set \\
			$V_{i^+}$, $V_{i^-} $ & set of pick-up nodes, set of drop-off nodes corresponding to request $i$ \\
			$A$ & arc set \\
			$c_a$ &  routing cost on arc $a$ \\
			$t_a$ &  travel time on arc $a$ \\
			$t_i$ &  travel time along the shortest path for request $i$  \\
			$\delta^{\text{in}}(v)$, $\delta^{\text{out}}(v)$  & incoming arcs, outgoing arcs of node $v$ \\
			$TW$ & fixed length of time window constructed from desired pick-up or desired drop-off time \\
			\bottomrule
		\end{tabularx}
		\caption{List of parameters.}
		\label{tab:paras}
	\end{threeparttable}
\end{table}
%
%
\begin{table}
	\centering
	\begin{threeparttable}
		\small
		\begin{tabularx}{0.8\textwidth}{cX}
			\toprule
			Variable & Description\\
			\midrule
			$p_i$ & binary variable indicating if user $i$ is transported or not \\
			$B_v$ & continuous variable indicating the start of service time at node $v$ \\
			$x_a$ & binary variable indicating if arc $a$ is used or not \\
			$d_i$ & regret of user $i$ w.r.t.\ $e_{i^-}$ \\ 
			$d_{\max}$ & maximum regret\\
			\bottomrule
		\end{tabularx}
		\caption{List of variables.}
		\label{tab:vars}
	\end{threeparttable}
\end{table}

Based on the event-based graph model, we are now ready to formulate our first MILP model for the DARP.

\subsection{Basic MILP for the  DARP}
\noindent
In this subsection, we propose an event-based mixed integer linear program for the DARP. First, we present a nonlinear mixed integer programming formulation, which is based on the event-based graph model presented in Section~\ref{subsec:event} above. In a second step, this model is transformed into an MILP by a reformulation of time window and ride time constraints,  i.e., constraints involving the variables $B_v$, $v\in V$, using a big-M method.

The DARP can be formulated as the following nonliner mixed integer program:
\begin{subequations}
	\begin{align}
		\min \;& \sum_{a\in A} c_a\,x_a  \label{costobj_first}\\[1ex]
		\text{\textit{s.\,t.}}\;
		&\sum_{a\in \delta^{\text{in}}(v)} \!\!\!x_a - \sum_{a\in \delta^{\text{out}}(v)}\!\!\! x_a = 0 \quad \forall v\in V, \label{flowconservation1} \\
		& \sum_{\substack{a\in\delta^{\text{in}}(v)\\ v \in V_{i^+}}}\!\!\! x_a = 1 \quad \forall i \in R, \label{passpickedup1} \\
		& \sum_{a \in \delta^{\text{out}}(\boldsymbol{0})} \! \! \! \!\! \! x_a \leq \vert K\vert, \label{numberofvehicles1} \\
		& B_w  \geq (B_v + s_{v_1} +  t_{(v,w)})\,x_{(v,w)} \quad \forall (v,w) \in A,\, v\neq \boldsymbol{0} \label{traveltime_nonlin} \\
		& B_w \geq e_0 + t_{(v,w)}\,x_{(v,w)} \quad \forall (\boldsymbol{0},w) \in A, \label{traveltime_nonlinb} \\
		& e_j \leq B_v \leq \ell_j \quad \forall j \in P \cup D \cup \{0\},\, v\in V_j \label{timewindows} \\
		& (B_w  - B_v -  s_{i^+})\!\! \sum_{a \in \delta^{\text{in}}(v)}\!\!\! x_a\! \! \sum_{a \in \delta^{\text{in}}(w)}\!\!\! x_a \leq L_i \quad \forall i \in R,\, v \in V_{i^+}, \, w \in V_{i^-}, \label{ridetime_nonlin} \\
		& x_a \in \{0,1\} \quad \forall a \in A, \\
		& B_v \geq 0 \quad \forall v \in V.
	\end{align}
\end{subequations}
The objective function \eqref{costobj_first} minimizes the total routing cost. While constraints \eqref{flowconservation1} are flow conservation constraints, it is ensured by constraints \eqref{passpickedup1} that each request $i\in R$ is accepted and that exactly one node of all nodes which contain the request's  pick-up location is reached by exactly one vehicle. Together with constraints \eqref{flowconservation1} and \eqref{passpickedup1}, the number of feasible dicycles in $G$ is bounded from above by the number of vehicles in constraints \eqref{numberofvehicles1}. Constraints \eqref{traveltime_nonlin} and \eqref{traveltime_nonlinb} define the difference in time needed to travel from one node to another. For all nodes in $V$ the start of service has to take place within the time window corresponding to the associated location of the node, which is handled by constraints \eqref{timewindows}. An upper bound on the ride time is ensured by constraints \eqref{ridetime_nonlin}. Note that we only impose a bound on the variables $B_w$, $B_v$, if both $v\in V_{i^+}$ and $w\in V_{i^-}$ are in fact the pick-up and drop-off nodes that are used to serve request $i$. Vehicle capacity, pairing and precedence constraints are ensured by the structure of the underlying network.

This formulation is nonlinear due to constraints \eqref{traveltime_nonlin} and \eqref{ridetime_nonlin}. In the following MILP these constraints are substitued by a linearized reformulation: 

\begin{model}\label{mod2}
	\begin{subequations}
		\begin{align}
			\min \;& \sum_{a\in A} c_a\,x_a  \label{costobj1}\\[1ex]
			\text{s.\,t.}\;  &\text{constraints}\ \eqref{flowconservation1} - \eqref{numberofvehicles1} \notag \\
			& B_w \geq B_v + s_{v_1} + t_{(v,w)} - \tilde{M}_{v,w}\,(1-x_{(v,w)}) \quad \forall (v,w) \in A,\, v \neq \boldsymbol{0}, \label{traveltime1} \\
			&B_w \geq e_0 +  t_{(v,w)}\,x_{(v,w)}  \quad \forall (\boldsymbol{0},w) \in A, \label{traveltime1b} \\
			& e_j \leq B_v \leq \ell_j \quad \forall j \in P \cup D \cup \{0\},\, v\in V_j \label{timewindows1} \\
			&B_w - B_v -  s_{i^+} \leq L_i + M_i\,\biggl(1-\sum_{a \in \delta^{\text{in}}(v)}x_a\!\! + 1-\sum_{a \in \delta^{\text{in}}(w)}x_a \biggr)  \quad \forall i \in R,\, v \in V_{i^+}, \, w \in V_{i^-}, \label{ridetime1} \\
			& x_a \in \{0,1\} \quad \forall a \in A, \\
			& B_v \geq 0 \quad \forall v \in V,
		\end{align}
	\end{subequations}
	where $M_i\geq \ell_{i^-} - e_{i^+} - L_i -  s_{i^+}$ and $\tilde{M}_{v,w}\geq \ell_{v_1} - e_{w_1} + s_{v_1} + t_{(v,w)}$ are sufficiently large constants. 
\end{model}
To include the option to deny requests in the DARP,  variables $p_i \in \lbrace 0, 1\rbrace$, $i \in R$ have to be added to Model~\ref{mod2} and constraints \eqref{passpickedup1} have to be changed to
\begin{equation}
	\sum_{\substack{a\in\delta^{\text{in}}(v)\\ v \in V_{i^+}}}\!\!\! x_a = p_i \qquad \forall i \in R. \label{passpick?}
\end{equation}
Hence, if a user is not picked-up (i.e., if $p_i=0$), then none of the nodes which contain his or her pick-up location are traversed by any vehicle. Note that in this case a reasonable objective function (see Section~\ref{sec:obj}) has to penalize the denial of user requests since otherwise an optimal solution is given by $p = 0$, $x=0$ and $B_v = e_{v_1}$ for all $v\in V$. In the computational experiments in Section~\ref{sec:experiments} we consider both cases, i.e., the scenario that all users have to be served and the scenario that some requests may be denied.

Assuming $q_i=1$ for all requests $i\in R$ and $n\geq Q+1$ the total number of variables in Model~\ref{mod2} can be bounded by $\mathcal{O}(n^{Q+1})$ with $\mathcal{O}(n^{2Q-1})$ constraints, of which $\mathcal{O}(n^{2Q-1})$ constraints are ride time constraints \eqref{ridetime1}. If $q_i \in \{2, \ldots, Q\}$ for some requests, the number of variables and constraints decreases. 

In each of the ride time constraints in Model~\ref{mod2}, the sums $\sum_{a \in \delta^{\text{in}}(v)}x_a $ and $\sum_{a \in \delta^{\text{in}}(w)}x_a$ are evaluated. This is computationally expensive, as will be demonstrated in the computational tests carried out in Section~\ref{sec:experiments}. By taking advantage of the relationship between the pick-up and drop-off time windows associated with request $i$, we show in the following how constraints \eqref{ridetime1} can be reformulated without using big-M constraints, resulting in a second MILP formulation referred to as Model~\ref{mod3}.

\subsection{Reformulation of time-related constraints}
The MILP model presented in this subsection differs from the previous model, Model~\ref{mod2}, in the formulation of time windows and ride time constraints. In Model~\ref{mod2}, the ride time constraints are modeled as big-M constraints that are used to deactivate the respective constraints for pick-up and drop-off nodes that are not contained in a vehicle's tour. By reformulating the time window constraints and using the relationship between earliest pick-up and latest drop-off times, the numerically unfavorable big-M constraints can be replaced by simpler constraints in the following model (Model~\ref{mod3}). This model is faster to solve which is verified by the numerical experiments presented in Section~\ref{sec:experiments}.

In this model, the ride time constraints, which ensure that a user does not spend more than $L_i$ minutes in the vehicle, are given by

\begin{equation}
	B_w - B_v -  s_{i^+} \leq L_i \qquad \forall i \in R,\,\forall v \in V_{i^+}, \,\forall w \in V_{i^-}. \label{ridetime_test}
\end{equation}

To show that these constraints, together with a reformulation of constraints \eqref{timewindows1}, reflect the modeling assumptions, we first observe that in general applications of  the DARP, as described for example in \cite{Cordeau2006}, users often formulate \emph{inbound} requests and \emph{outbound} requests. In the first case, users specify a desired departure time from the origin, while in the case of an outbound request, users specify a desired arrival time at the destination. In both cases a time window of a fixed length $TW$ is constructed from the desired time, so that we end up with a pick-up time window of length $TW = \ell_{i^+}-e_{i^+}$ for an inbound request and a drop-off time window of length $TW=\ell_{i^-}-e_{i^-}$ for an outbound request. Now the remaining time window is constructed as follows (based on \cite{Cordeau2006}): For an inbound request the drop-off time window is given by the bounds
\begin{equation}
	e_{i^-} = e_{i^+} +  s_{i^+} + t_{i} \qquad \text{and} \qquad\ell_{i^-} = \ell_{i^+} +  s_{i^+} + L_i, \label{inbound} 
\end{equation}
where $t_i$ denotes the direct travel time $t_{(v,w)}$ from a node $v\in V_{i^+}$ to a node $w\in V_{i^-}$ , i.e.\ $t_i = t_{(v,w)}$ with $v = (i^+,0,\ldots,0)$ and $w = (i^-,0,\ldots,0)$. Similarly, for an outbound request the pick-up time window is defined by
\begin{equation}
	e_{i^+} = e_{i^-} - L_i -  s_{i^+} \qquad \text{and} \qquad \ell_{i^+} = \ell_{i^-} - t_i -  s_{i^+}. \label{outbound}
\end{equation}
Secondly, we define the notion of an \emph{active} node $v\in V$: We say that a node $v\in V$ is \emph{active} if at least one of its incoming arcs is part of a dicycle flow,  i.e., if
\begin{equation*}
	\sum_{a\in \delta^{\text{in}}(v)}x_a = 1.
\end{equation*}
Otherwise, we call $v$ \emph{inactive}. Note that due to constraints \eqref{passpickedup1}, for each request $i\in R$ we have exactly one associated active pick-up node and one associated active drop-off node. In case we include the option of denying user requests, i.e.\ we use constraints \eqref{passpick?} instead of constraints \eqref{passpickedup1} and add variables $p_i \in \{0,1\}$, $i\in R$ to the MILP, each request either has exactly one active pick-up and drop-off node (in this case we have $p_i=1$), or no associated active node at all (this is the case when $p_i=0$). 

Now, if both $v$ and $w$ in constraints \eqref{ridetime_test} are active nodes, inequalities \eqref{ridetime_test} and \eqref{ridetime1} coincide. In case both $v$ and $w$ are inactive, the values $B_v$, $B_w$ can be ignored in an interpretation of an optimal solution, as $v$ and $w$ are not contained in any of the vehicle tours. Hence, the critical two cases are the cases where one of the nodes is active and the other node is inactive. Let $v^{\text{off}}$ and $v^{\text{on}}$ denote the inactive and the active node from the set $\{v,w\}$, respectively. Then, we do not want that $B_{v^{\text{off}}}$  influences the value of  $B_{v^{\text{on}}}$  in constraints \eqref{ridetime_test}:
\begin{description}
	\item[Case 1: $v$ is active, $w$ is inactive.]  Resolving \eqref{ridetime_test} for $B_v$ we obtain $B_v \geq B_w - L_i -  s_{i^+}$. Now, we do not want to impose any additional constraints on $B_v$. Thus, we demand $B_w  - L_i -  s_{i^+} \leq e_{i^+}$. Accordingly, it has to hold that 
	\begin{equation*}
		B_w \leq e_{i^+} + L_i +  s_{i^+}.
	\end{equation*}
	Recall that here we assume $w$ to be inactive. If $w$ is active, the weaker constraint $B_w \leq \ell_{i^-}$ needs to hold. Putting these restrictions together for an inbound request, we get
	\begin{align}
		B_w &\leq e_{i^+} + L_i +  s_{i^+} + (\ell_{i^-}- (e_{i^+} + L_i +  s_{i^+}))\sum_{a \in \delta^{\text{in}}(w)} x_a \notag \\
		& =  e_{i^+} + L_i +  s_{i^+} + (\ell_{i^+} + L_i +  s_{i^+} - (e_{i^+} + L_i +  s_{i^+}))\sum_{a \in \delta^{\text{in}}(w)} x_a \notag  \\
		&=  e_{i^+} + L_i +  s_{i^+} + (\ell_{i^+}  - e_{i^+} )\sum_{a \in \delta^{\text{in}}(w)} x_a \notag\\
		&= e_{i^+} + L_i +  s_{i^+} + TW\sum_{a \in \delta^{\text{in}}(w)} x_a, \label{B_w_bound1}
	\end{align} 
	using the reformulation of $\ell_{i^-}$ from equations \eqref{inbound}. In the same manner, we use equations \eqref{outbound} to substitute $e_{i^+} = e_{i^-} - L_i -  s_{i^+}$ and obtain that
	\begin{equation}
		B_w \leq e_{i^+} + L_i +  s_{i^+} + (\ell_{i^-}  - e_{i^-} )\sum_{a \in \delta^{\text{in}}(w)} x_a  = e_{i^+} + L_i +  s_{i^+} + TW\sum_{a \in \delta^{\text{in}}(w)} x_a\label{B_w_bound2}
	\end{equation}
	has to hold for an outbound request. 
	Hence the two formulations \eqref{B_w_bound1} and \eqref{B_w_bound2} coincide. 
	
	\item[Case 2: $\boldmath{v}$ is inactive, $\boldmath{w}$ is active.]   Resolving \eqref{ridetime_test} for $w$ we obtain $B_w \leq L_i + B_v +  s_{i^+}$. In order for the latter inequality to be redundant for $B_w$, we demand that $L_i + B_v + s_{i^+} \geq \ell_{i^-}$. It follows that $B_v \geq \ell_{i^-} - L_i -  s_{i^+}$ has to hold. For an inbound request, using the equivalence from \eqref{inbound}, this can be resolved to 
	\begin{equation*}
		B_v \geq (\ell_{i^+} +  L_i + s_{i^+}) - L_i -  s_{i^+} = \ell_{i^+}. 
	\end{equation*}
	Recall that we assumed that $v$ is inactive. In case $v$ is active, the less tighter constraint $B_v \geq e_{i^+}$ has to hold, so that we arrive at 
	\begin{equation*}
		B_v\geq e_{i^+} +  (\ell_{i^+} - e_{i^+}) \,\biggl(1-\sum_{a \in \delta^{\text{in}}(v)} x_a\biggr) =e_{i^+} +  TW \,\biggl(1-\sum_{a \in \delta^{\text{in}}(v)} x_a\biggr). 
	\end{equation*}
	In a similar fashion, we obtain
	\begin{equation*}
		B_v \geq e_{i^+} +  (\ell_{i^-} - e_{i^-}) \,\biggl(1-\sum_{a \in \delta^{\text{in}}(v)} x_a\biggr) =  e_{i^+} + TW \,\biggl(1-\sum_{a \in \delta^{\text{in}}(v)} x_a\biggr) 
	\end{equation*}
	for an outbound request. 
	We conclude that the two lower bounds on $B_v$ coincide.
\end{description}
As desired, by reformulating the constraints \eqref{timewindows1} on the variables $B_v$, $v\in V$, we obtain a simpler version of the ride time constraints \eqref{ridetime1}. We put these results together in a second MILP formulation of the DARP.
\begin{model}\label{mod3}
	\begin{subequations}
		\begin{align}
			\min \;& \sum_{a\in A} c_a\,x_a  \label{costobj}\\[1ex]
			\text{s.\,t.}\;  &\text{constraints}\ \eqref{flowconservation1} - \eqref{numberofvehicles1} \notag \\
			& B_w \geq B_v + s_{v_1} + t_{(v,w)} - \tilde{M}_{v,w}\,(1-x{(v,w)}) \quad \forall (v,w) \in A,\, v\neq \boldsymbol{0} \label{traveltime} \\
			&B_w \geq e_0 +  t_{(v,w)}\,x_{(v,w)}  \quad \forall (\boldsymbol{0},w) \in A, \label{traveltimeb} \\
			& e_0 \leq B_{\boldsymbol{0}} \leq \ell_0,\label{timewindows_0}\\
			& e_{i^+} + TW\,\biggl(1- \sum_{a \in \delta^{\text{in}}(v)}\!\!x_a\biggr) \leq B_v \leq \ell_{i^+} \quad \forall i \in R,\, v\in V_{i^+}, \label{timewindows_pickup} \\
			& e_{i^-} \leq B_v \leq e_{i^+} + L_i +  s_{i^+} + TW \!\sum_{a \in \delta^{\text{in}}(v)}\!\!x_a \quad \forall i\in R,\,  v \in V_{i^-}, \label{timewindows_dropoff} \\
			& B_w  - B_v -  s_{i^+}\leq L_i \quad \forall i \in R,\, v \in V_{i^+}, \, w \in V_{i^-}, \label{ridetime} \\
			& x_a \in \{0,1\} \quad \forall a \in A, \\
			& B_v \geq 0 \quad \forall v \in V. \label{lastconstr}
		\end{align}
	\end{subequations}
\end{model}
We substitute ride time constraints \eqref{ridetime1} for a simpler version \eqref{ridetime}. In return we use a more complex version of the time window constraints \eqref{timewindows_0} -- \eqref{timewindows_dropoff} instead of the short form \eqref{timewindows1}. 

Similar to Model~\ref{mod2}, variables $p_i \in \{0,1\}$, $i\in R$ may be added to Model~\ref{mod3} and constraints \eqref{passpickedup1} may be substituted by constraints \eqref{passpick?} to include the option of denying user requests. In this case, the objective function should be modified to contain a term penalizing the denial of requests.

As there are $\mathcal{O}(n^{2Q-1})$ ride time constraints, for each of which two sums $\sum_{a \in \delta^{\text{in}}(v)} x_a$ have to be evaluated in the longer version \eqref{ridetime1}, but only $\mathcal{O}(n^Q)$ time window constraints, for each of which one sum of the above form has to be evaluated in the longer version \eqref{timewindows_0} - \eqref{timewindows_dropoff}, we obtain a more efficient new MILP formulation. Similar to Model~\ref{mod2}, for $n\geq Q+1$ there are at most $\mathcal{O}(n^{Q+1})$ variables and at most $\mathcal{O}(n^{2Q-1})$ constraints, of which $\mathcal{O}(n^{2Q-1})$ constraints are ride time constraints \eqref{ridetime}.

\subsection{Objective functions}\label{sec:obj}
\noindent
In most of the research on the DARP only one objective is used, which is often the minimization of total routing costs. An excellent overview is given by \cite{Ho2018}. Other popular objectives are, for example, the minimization of total route duration, number of vehicles used, users' waiting time, drivers' working hours, or deviation from the desired pick-up and drop-off times.

In this paper, we focus on three prevalent (and possibly conflicting) criteria, namely the total routing cost, the total number of unanswered requests and the total users' regret or the maximum regret, and combine these three criteria into weighted sum objective functions. The first and probably most important criterion is the total routing cost, which can be computed as 
\begin{equation}
	f_c(x) \coloneqq \sum_{a \in A} c_a \, x_a. \label{f_c} 
\end{equation}
We refer to $f_c$ as \emph{cost-objective}. The second objective function,
\begin{equation*}
	f_{n} (p)\coloneqq n - \sum_{i\in R} p_i,
\end{equation*}
measures the total number of declined requests. The next optimization criterion relates to customer satisfaction: We measure the response time to a service request by assessing a user's regret, which aims at penalizing overly long travel times as well as possibly delayed pick-up times.

Let the variable $d_i \geq 0$, $i \in R$ measure the difference in time compared to a user's earliest possible arrival time. We refer to $d_i$ as a user's \emph{regret}. Moreover, let the variabe $d_{\max} \geq 0$ measure the \emph{maximum regret}. By introducing constraints 
\begin{align}
	d_i  &\geq  B_v - e_{i^-} \qquad\forall i\in R,\, \forall v \in V_{i^-} , 	\label{regret}\\
	d_{\max}  & \geq d_i \qquad \forall i \in R, \label{maxregret}
\end{align}
we can now minimize the total or average regret, or the maximum regret, respectively. The total regret is thus given by the \emph{regret-objective} 
\begin{equation}
	f_r(d) \coloneqq \sum_{i \in R} d_i, \label{f_r}
\end{equation}
while the \emph{maximum-regret-objective} is given by 
\begin{equation}
	f_{r_{\max} }(d_{\max})\coloneqq d_{\max}. \label{f_rmax}
\end{equation}

The discussion above highlights the fact that we have to consider  different and generally conflicting objective functions that are relevant when solving the DARP. While the cost objective $f_c$ aims at minimizing total travel cost and thus takes the perspective of the service provider, the quality of service which rather takes a user's perspective is better captured by objective functions like $f_n$, $f_e$ and $f_{e_{\max}}$. 

We approach this technically bi- or multi-objective problem by using a weighted sum approach with fixed weights, i.e., by combining the relevant objective functions into one weighted sum objective.
When using the total regret as quality criterion, we obtain
\begin{equation}
	f_{cr} (x, d)\coloneqq \sum_{a \in A} c_a \, x_a  + \alpha\sum_{i \in R} d_i \label{f_cr}
\end{equation}
which will be referred to as \emph{cost-regret-objective}, and when using the maximum regret as quality criterion, we get
\begin{equation}
	f_{cr_{\max}} (x, d_{\max})\coloneqq \sum_{a \in A}c_a \, x_a + \beta d_{\max} \label{f_crmax}
\end{equation}
which will be referred to as \emph{cost-max-regret-objective}. The parameters $\alpha>0$ and $\beta>0$ are weighting parameters that can be selected according to the decision maker's preferences. We refer to \cite{ehrgott05multicriteria} for a general introduction into the field of multi-objective optimization.

Last but not least, we consider a form of the DARP in which it is allowed to deny certain user requests. 
This is accomplished by substituting constraints \eqref{passpickedup1} in Model~\ref{mod2}  or Model~\ref{mod3} , respectively, by constraints \eqref{passpick?} and adding variables $p_i \in \{0,1\}$, $i\in R$ to Model~\ref{mod2} or \ref{mod3}. In this case, the number of accepted requests has to be maximized or, equivalently, the number of unanswered requests has to be minimized. At the same time, routing costs and regret should be as small as possible. The optimization of these opposing criteria is reflected by the \emph{request-cost-regret-objective} given by
\begin{equation}
	f_{rcr}(x, d, p) \coloneqq \sum_{a \in A} c_a \, x_a  + \alpha\sum_{i \in R} d_i + \gamma \left(n - \sum_{i \in R} p_i \right), \label{f_rcr}
\end{equation}
where $\gamma>0$ is an additional weighting parameter. While the third part of the objective refers to the number of unanswered requests and is equal to $\gamma\, n $ at maximum (i.e., if no requests are accepted), the values of the total routing costs and of the total regret strongly depend on the underlying network and request data. Note that meaningful choices of the weighting parameters have to reflect this in order to avoid situations where one part of the objective overrides the others. This will be discussed in Section~\ref{sec:experiments}. 

\section{Numerical Experiments}\label{sec:experiments}
\noindent
This section is divided into two parts. In the first subsection, we compare the MILP formulations  Model~\ref{mod2} and \ref{mod3} to a standard compact 3-index formulation of  the DARP from \cite{Cordeau2006} and to a compact 2-index formulation for the PDPTW by \cite{Furtado2017}, which can be adapted to a MILP formulation for the DARP by using additional constraints. The computational performance is evaluated on a set of benchmark test instances.
In the second part, we substitute the cost objective function in Model~\ref{mod3} with different objective functions as introduced in Section~\ref{sec:obj} and analyze the effect with respect to economic efficiency and customer satisfaction. For this purpose, a set of $60$ artificial instances from the city of Wuppertal is generated. The computations are carried out on an Intel Core i7-8700 CPU, 3.20\,GHz, 32\,GB memory using CPLEX 12.10. The MILPs are programmed using C++. The time limit for the solution in all tests is set to $7200$ seconds. 
In the following, the computational results are discussed in detail. We use the abbreviations listed in Table~\ref{tab:abrev}.

\begin{table}
	\small
	\begin{tabularx}{\textwidth}{lX}\toprule
		Inst. & Name of instance \\
		$n$ & Number of users in the respective instances. \\
		Gap & Gap obtained by CPLEX solution (in percent). \\
		CPU & Computational time in seconds. \\
		Obj.v. & Objective value. \\
		$f_c$ & Total routing costs. \\
		$f_r$ & Total regret. \\
		$f_{r_{\max}}$ & Maximum regret. \\
		a.r. & Answered requests (in percent). \\	
		$\Delta f_c$ & Change in total routing costs (in percent). \\
		$\Delta f_r$ & Change in total regret (in percent).\\
		$\Delta f_{r_{\max}}$ & Change in maximum regret (in percent). \\
		$\Delta$ a.r.  & Change in answered requests (in percent). \\\bottomrule
	\end{tabularx}
	\caption{List of abbreviations.\label{tab:abrev}}
\end{table}

When computing the average CPU times over several runs, and an instance was not solved to optimality within the time limit, then a CPU 
time of 7200 seconds is assumed.

\subsection{Benchmark data}
\noindent
We compare our MILP models to a standard compact 3-index formulation of the DARP from the literature. For this purpose, we use the basic mathematical model of the DARP introduced by \cite{Cordeau2006}. Note that a tighter 2-index formulation is given by \cite{Ropke2007}. However, this comes at the price of an exponentially growing number of constraints. It is thus better suited for a solution within a B\&C framework and we did not include it in our comparison. Moreover, we compare our model to a compact 2-index formulation for the PDPTW, introduced by \cite{Furtado2017}, which can be adapted to a MILP formulation of the DARP by adding constraints on the passengers' maximum ride time, the maximum duration of service and the maximum number of vehicle tours.

Both models, introduced by \cite{Cordeau2006} and \cite{Furtado2017}, referred to as C-DARP and F-DARP in the following, are based on a complete directed graph. The node set comprises all pick-up and drop-off locations and two additional nodes $0$ and $2n+1$ for the depot. Thus, the node set is equal to the set $P\cup D \cup \{0, 2n+1\}$. For C-DARP, we use the MILP formulation with a reduced number of variables and constraints, described in \cite{Cordeau2006}, which includes aggregated variables $B_j$ (modelling the beginning of service time) and $Q_j$ (modelling the vehicle load)  at every node $j$ except the origin and destination depot.  The objective is to minimize the total routing costs. We do not add any additional valid inequalities described in \cite{Cordeau2006} to any of the MILP formulations in this comparison. We use the following trivial variable fixings: $x_{i,i} = 0$, $x_{i,0} = 0$ and $x_{2n+1,i} = 0$ for all $i\in P\cup D$, $x_{0,n+i} =0$, $x_{i,2n+1} = 0$ and $x_{n+i,i} = 0$ for all $i\in P$. Besides that, we use time window, capacity, and ride time constraints to fix the following variables: Let $i,j\in P$, then $x_{i,j}=0$ if $f^1_{j,i} = f^2_{j,i} = 0$, $x_{n+i,n+j} = 0$ if $f^1_{j,i} = f^2_{i,j} = 0$, and $x_{i,n+j} = 0$ if $f^1_{i,j} = 0$.

The test instances are the two sets of benchmark instances\footnote{The instances are available at \url{http://neumann.hec.ca/chairedistributique/data/darp/branch-and-cut/}.} set $a$ and set $b$ created in \cite{Cordeau2006}. 
The characteristics of the instances are summarized in Table~\ref{tab:charbm}. 
In all test instances we tighten the remaining time windows, i.e.\ the time windows not given by the pick-up time of inbound requests or by the drop-off time of outbound requests, respectively, as described in \cite{Cordeau2006}: The bounds of the missing time windows can be calculated according to equations \eqref{inbound} and \eqref{outbound} stated earlier in Section~\ref{sec:model}.

\begin{table}
	\centering
	\begin{threeparttable}
		\footnotesize
		\begin{tabular}{ccccccp{4ex}cccccc}
			\cmidrule(lr){1-6}   \cmidrule(lr){8-13}
			Instance & $Q$ & $n$ & $\vert K\vert$ & $L_i$ & $T$ && Instance & $Q$ & $n$ & $\vert K\vert$ & $L_i$ & $T$ \\ 
			\cmidrule(lr){1-6}   \cmidrule(lr){8-13}
			a2-16& 3 & 16 & 2 &	30 & 480 && b2-16 & 6 & 16 & 2 & 45 & 480 \\
			a2-20& 3 & 20 & 2 &	30 & 600 && b2-20 & 6 & 20 & 2 & 45 & 600 \\
			a2-24& 3 & 24 & 2 & 30 & 720 && b2-24 & 6 & 24 & 2 & 45 & 720 \\
			a3-18& 3 & 18 & 3 & 30 & 360 && b3-18 & 6 & 18 & 3 & 45 & 360 \\
			a3-24& 3 & 24 & 3 & 30 & 480 && b3-24 & 6 & 24 & 3 & 45 & 480 \\
			a3-30& 3 & 30 & 3 & 30 & 600 && b3-30 & 6 & 30 & 3 & 45 & 600 \\
			a3-36& 3 & 36 & 3 & 30 & 720 && b3-36 & 6 & 36 & 3 & 45 & 720 \\
			a4-16& 3 & 16 & 4 & 30 & 240 && b4-16 & 6 & 16 & 4 & 45 & 240 \\
			a4-24& 3 & 24 & 4 & 30 & 360 && b4-24 & 6 & 24 & 4 & 45 & 360 \\
			a4-32& 3 & 32 & 4 & 30 & 480 && b4-32 & 6 & 32 & 4 & 45 & 480 \\
			a4-40& 3 & 40 & 4 & 30 & 600 && b4-40 & 6 & 40 & 4 & 45 & 600 \\
			a4-48& 3 & 48 & 4 & 30 & 720 && b4-48 & 6 & 48 & 4 & 45 & 720 \\
			a5-40& 3 & 40 & 5 & 30 & 480 && b5-40 & 6 & 40 & 5 & 45 & 480 \\
			a5-50& 3 & 50 & 5 & 30 & 600 && b5-50 & 6 & 50 & 5 & 45 & 600 \\
			a5-60& 3 & 60 & 5 & 30 & 720 && b5-60 & 6 & 60 & 5 & 45 & 720 \\
			a6-48& 3 & 48 & 6 & 30 & 480 && b6-48 & 6 & 48 & 6 & 45 & 480 \\
			a6-60& 3 & 60 & 6 & 30 & 600 && b6-60 & 6 & 60 & 6 & 45 & 600 \\
			a6-72& 3 & 72 & 6 & 30 & 720 && b6-72 & 6 & 72 & 6 & 45 & 720 \\
			a7-56& 3 & 56 & 7 & 30 & 480 && b7-56 & 6 & 56 & 7 & 45 & 480 \\
			a7-70& 3 & 70 & 7 & 30 & 600 && b7-70 & 6 & 70 & 7 & 45 & 600 \\
			a7-84& 3 & 84 & 7 & 30 & 720 && b7-84 & 6 & 84 & 7 & 45 & 720 \\
			a8-64& 3 & 64 & 8 & 30 & 480 && b8-64 & 6 & 64 & 8 & 45 & 480 \\
			a8-80& 3 & 80 & 8 & 30 & 600 && b8-80 & 6 & 80 & 8 & 45 & 600 \\
			a8-96& 3 & 96 & 8 & 30 & 720 && b8-96 & 6 & 96 & 8 & 45 & 720 \\
			\cmidrule(lr){1-6}   \cmidrule(lr){8-13}
		\end{tabular}
		\caption{Characteristics of the benchmark test instances.}
		\label{tab:charbm}
	\end{threeparttable}
\end{table}
\begin{table}
	\centering
	\begin{threeparttable}
		\footnotesize
		\begin{tabular}{cccccccccccc}
			\toprule
			& \multicolumn{3}{c}{C-DARP} & \multicolumn{3}{c}{F-DARP}& \multicolumn{2}{c}{Model~\ref{mod2}} & \multicolumn{2}{c}{Model~\ref{mod3}} \\
			\cmidrule(lr){2-4}\cmidrule(lr){5-7}\cmidrule(lr){8-9}\cmidrule(lr){10-11}\addlinespace[0.1em]
			Inst. & Obj.v.& Gap & \rev{CPU}  & Obj.v.& Gap & \rev{CPU}  & Obj.v.  & \rev{CPU} &  Obj.v.  & 	\rev{CPU} \\
			\midrule
			a2-16 & 294.3 &          &0.08& 294.3 &      &0.04 &  294.3&  0.04& 294.3&0.10\\
			a2-20 & 344.9 &          &0.23& 344.9 &      &0.12  &  344.9&  0.02& 344.9&0.02\\
			a2-24 & 431.1 &          &0.68& 431.1 &      &0.28  &  431.1&  0.06& 431.1&0.05\\
			a3-18 & 300.5 &          &0.68& 300.5 &      &0.24  &  300.5&  0.05& 300.5&0.04\\
			a3-24 & 344.9 &          &5.53& 344.9 &      &1.93  &  344.9&  0.10& 344.9&0.10\\
			a3-30 & 494.8 &          &20.14& 494.8 &      &297  &  494.8&  0.07& 494.8&0.04\\
			a3-36 & 583.2 &          &28.70& 583.2 &      &6.65  &  583.2&  0.09& 583.2&0.11\\
			a4-16 & 282.7 &          &5.61& 282.7 &      &2.85   &  282.7&  0.05& 282.7&0.06\\
			a4-24 & 375.0 &          &11.13& 375.0 &      &24.95 &  375.0&  0.07& 375.0&0.04\\
			a4-32 & 485.5 &          &591& 485.5 & 8  &7200  &  485.5&  0.12& 485.5&0.15\\
			a4-40 & 557.7 &          &1903& 557.7 & 4  & 7200 &  557.7&  0.48& 557.7&0.28\\
			a4-48 & 669.4 & 6      &7200& 680.0 & 10 & 7200  &  668.8&  0.67& 668.8&0.50\\
			a5-40 & 498.4 & 2      &7200& 498.4 & 11 & 7200 &  498.4&  0.23& 498.4&0.23\\
			a5-50 & 696.3 & 17    &7200& N/A    & N/A &7200 &  686.6&  4.02& 686.6&2.47\\
			a5-60 & 828.3 & 14    &7200& N/A    & N/A &7200 &  808.4&  1.41& 808.4&1.03\\
			a6-48 & 604.1 & 18    &7200& N/A    & N/A &7200 &  604.1&  1.33& 604.1&0.89\\
			a6-60 & 874.8 & 16    &7200& N/A    & N/A &7200 &  819.3&  6.17& 819.3&8.83\\
			a6-72 &  N/A    & N/A &7200& N/A    & N/A &7200 &  916.1&  17.88& 916.1&13.87\\
			a7-56 & 764.5 & 19    &7200& N/A    & N/A &7200 &  724.0&  1.92& 724.0&1.24\\
			a7-70 &  N/A   & N/A &7200& N/A & N/A     &7200 &  875.7&  4.77& 875.7&6.93\\
			a7-84 & N/A    & N/A &7200& N/A & N/A      &7200 &  1033.3&  35.70& 1033.3&5.89\\
			a8-64 & N/A    & N/A &7200& N/A & N/A     &7200 &  747.5&  15.31& 747.5&7.19\\
			a8-80 & N/A    & N/A &7200& N/A & N/A     &7200 &  945.8& 51.19 & 945.8&25.11\\
			a8-96 & N/A    & N/A&7200& N/A & N/A      &7200 &  1229.7&  3593& 1229.7&461\\
			\bottomrule 
		\end{tabular}
		\begin{tablenotes}
			\small
			\item[N/A ]Not applicable. No integer solution found within the time limit.
		\end{tablenotes}
		\caption{Solution values and computing times for the benchmark test instances 'a'.}
		\label{tab:ainstances}
	\end{threeparttable}
\end{table}
\begin{table}
	\centering
	\begin{threeparttable}
		\footnotesize
		\begin{tabular}{cccccccccccc}
			\toprule
			& \multicolumn{3}{c}{C-DARP} & \multicolumn{3}{c}{F-DARP}& \multicolumn{2}{c}{Model~\ref{mod2}} & \multicolumn{2}{c}{Model~\ref{mod3}} \\
			\cmidrule(lr){2-4}\cmidrule(lr){5-7}\cmidrule(lr){8-9}\cmidrule(lr){10-11}\addlinespace[0.1em]
			Inst. & Obj.v. & Gap & \rev{CPU} & Obj.v. & Gap & \rev{CPU} & Obj.v.  & \rev{CPU}  &  Obj.v.  & \rev{CPU} \\
			\midrule
			b2-16 &309.4&            & 0.13 &309.4 &         & 0.13& 309.4& 0.03&309.4&0.03\\
			b2-20 &332.7&            & 0.06 & 332.7&         & 0.02& 332.7& 0.02&332.7&0.01\\
			b2-24 &444.7&            & 0.91 & 444.7&         & 0.81& 444.7& 0.06&444.7&0.05\\
			b3-18 &301.6&            & 0.83 & 301.6&         & 0.73& 301.6& 0.05&301.6&0.04\\
			b3-24 &394.5&            & 5.28 & 394.5&         & 0.74& 394.5& 0.11&394.5&0.11\\
			b3-30 &531.4&            & 1.68 & 531.4&         & 0.22& 531.4& 0.07&531.4&0.06\\
			b3-36 &603.8&            & 6.41 & 603.8&         & 2.17& 603.8& 0.07&603.8&0.09\\
			b4-16 &297.0&            & 0.23 & 297.0&         & 0.03& 297.0& 0.04&297.0&0.02\\
			b4-24 &371.4&            & 1.86 & 371.4&         & 0.66& 371.4& 0.06&371.4&0.06\\
			b4-32 &494.9&            & 3.01 & 494.9&         & 0.57& 494.9& 0.04&494.9&0.04\\
			b4-40 &656.6&            & 51.15 & 656.6&      & 20.04& 656.6& 0.13&656.6&0.14 \\
			b4-48 &673.8& 2        & 7200 & 680.7& 6& 7200& 673.8& 1.14&673.8&0.94\\
			b5-40 &613.7& 10     & 7200& 619.1& 6&7200& 613.7& 0.20&613.7&0.21\\
			b5-50 &763.0& 7       & 7200 & 768.2& 8&7200& 761.4& 0.54&761.4&0.45\\
			b5-60 &917.6& 14     & 7200 & N/A& N/A&7200& 902.0& 4.01&902.0&0.92\\
			b6-48 &714.8&           & 2700 &715.9& 2&7200& 714.8& 0.27&714.8&0.26\\
			b6-60 &893.4& 11     &7200 &N/A&N/A &7200 &860.0& 0.50&860.0&0.43\\
			b6-72 &1083.9& 22  & 7200 &N/A&N/A&7200& 978.5& 3.48&978.5&5.11\\
			b7-56 &850.7& 14     & 7200 &867.9&14&7200& 824.0& 6.99&824.0&3.62\\
			b7-70 &919.7& 10     & 7200 &N/A&N/A &7200& 912.6& 2.54&912.6&2.42\\
			b7-84 & N/A  & N/A  & 7200 &N/A&N/A&7200& 1203.4& 2.93 &1203.4&2.72\\
			b8-64 & N/A  & N/A  & 7200 &N/A&N/A&7200& 839.9&2.03&839.9&2.16\\
			b8-80 &1433.0& 37  & 7200 &N/A&N/A&7200& 1036.4& 1.34&1036.4&1.54\\
			b8-96 & N/A  & N/A  & 7200 &N/A&N/A&7200& 1185.6& 27.42 &1185.6&24.26\\
			\bottomrule 
		\end{tabular}
		\begin{tablenotes}
			\small
			\item[N/A ]Not applicable. No integer solution found within the time limit.
		\end{tablenotes}
		\caption{Solution values and computing times for the benchmark test instances 'b'.}
		\label{tab:binstances}
	\end{threeparttable}
\end{table}

A summary of the computational results can be found in Tables~\ref{tab:ainstances} and \ref{tab:binstances}. For each of the considered models C-DARP, F-DARP, Model~\ref{mod2} and \ref{mod3}, the objective value (Obj.v.) of the cost objective and the computational time in seconds (CPU) are reported. 
For C-DARP and F-DARP we also report the relative gap (Gap) as some of the instances have not been solved to optimality within the time limit of 7200 seconds. By taking a closer look at Tables~\ref{tab:ainstances} and \ref{tab:binstances}, it becomes evident that Model~\ref{mod3} outperforms Model~\ref{mod2} in a majority of the instances, especially in the larger a-instances (in terms of number of users and vehicles used). Moreover, one can see clearly from the results that Models~\ref{mod2}  and \ref{mod3} yield a more efficient formulation than C-DARP and F-DARP: Starting from instance a4-48 (C-DARP) and instance a4-32 (F-DARP), the a-instances could not be solved to optimality within two hours (or even no integer solution was found). The computational time needed to solve these instances with Model~\ref{mod3} range from less than one second to 461 seconds and from less than one second to about 1 hour for Model~\ref{mod2}. In general, with Models~\ref{mod2} and \ref{mod3} the b-instances are easier to solve than the a-instances, as a large majority of the instances have been solved in less than one second. This reflects the fact that the size of the MILPs decreases tremendously when users request more than one seat, which can be traced back to the fact that the size of the event-based graph decreases in this case. With C-DARP and F-DARP, CPLEX was not able to solve instances b4-48 to b8-96 to optimality within the time limit of two hours, sometimes no integer solution was found. This clearly emphasizes the computational efficiency of the event-based graph as the underlying structure of an MILP formulation for the DARP.

\subsection{Artificial data -- City of Wuppertal}
\noindent
Ridepooling services are usually operated by taxis or mini-buses, whose passenger capacity is often equal to three or six. Thus, in the artificially created instances we consider the case that $Q\in\{3,6\}$. In the case that $Q=3$, we assume that each user requests one seat each, while for $Q=6$ the number of requested seats per user is chosen randomly from $\{1,\ldots, 6\}$. Moreover, the service time $s_j$ associated with location $j$ is set to be equal to the number of requested seats $q_j$. This is in accordance with the benchmark test instances for the DARP created in \cite{Cordeau2006}. 
The instance size is determined by the number of users. For both $Q=3$ and $Q=6$ and for each number of users $n\in \lbrace  20, 30, 40, 60, 80, 100\rbrace$ we generate a set of $5$ instances with $n$ users each. 
We denote the instances by $Q3.n.m$ and $Q6.n.m$, indicating the $m$-th instance with $n$ users, $m\in\lbrace1, \ldots,5\rbrace$. 
Pick-up and drop-off locations are chosen randomly from a list of streets in Wuppertal, Germany. The taxi depot is chosen to be located next to the main train station in Wuppertal. The cost $c_a$ for an arc $a$ in the event-based network is computed as the length of a shortest path from its tail to its head in an OpenStreetMap network corresponding to the city of Wuppertal. The shortest path is calculated based on OpenStreetMap data using the shortest path method of the Python package NetworkX. Due to slowly moving traffic within the city center, the average travel speed is set to 15~km/h, so that the travel time in minutes is equal to $t_a = 4 \,c_a$. The maximum duration of service $T$ is set to $240$ minutes. Earliest pick-up times are chosen randomly from five-minute intervals within the next 5--205 minutes (we consider inbound requests only) and the time window length is chosen to be equal to 15~minutes, as we assume that users of ridepooling services in a city, where public transport operates at high frequencies, want to be picked-up without long waiting times. A user $i$'s maximum ride time $L_i$ is set to $1.5$ times the ride time for a direct ride from the pick-up to the drop-off location. A summary of the artificial instances' remaining characteristics can be found in Table~\ref{tab:characteristics}.

\begin{table}
	\centering
	\begin{threeparttable}
		\footnotesize
		\begin{tabular}{cccccccc}
			\toprule
			Instances & $Q$ & $n$ & $\vert K\vert$ & Instances & $Q$ & $n$ & $\vert K\vert$ \\ 
			\midrule
			Q3.20.1-5& 3 & 20 &  4-5 & Q6.20.1-5& 6 & 20 & 5-6 \\
			Q3.30.1-5& 3 & 30 & 5-6  & Q6.30.1-5& 6 & 30 & 6-7 \\
			Q3.40.1-5& 3 & 40 & 6-7 & Q6.40.1-5& 6 & 40 &  8-9\\
			Q3.60.1-5& 3 & 60 & 9-10 & Q6.60.1-5& 6 & 60 & 11-12 \\
			Q3.80.1-5& 3 & 80 & 10-11 & Q6.80.1-5& 6 & 80 & 13-14 \\
			Q3.100.1-5& 3 & 100 & 13-14 & Q6.100.1-5& 6 & 100 & 16-17 \\
			\bottomrule
		\end{tabular}
		\caption{Characteristics of the Wuppertal artificial test instances.}
		\label{tab:characteristics}
	\end{threeparttable}
\end{table}

In general, the artificial instances are more complex than the benchmark instances: there is a smaller planning horizon (240--720 minutes for the benchmark instances and 240 minutes for the artificial instances) during which the same number of user requests have to be served, i.e.\ the ratio of user requests per time is higher. This is also reflected by the number of required vehicles. While there are only two to eight vehicles in the benchmark test instances to serve between $16$ and $96$ user requests, there are four to $17$ vehicles required in the artificial instances.

It has been shown in the previous subsection that Model~\ref{mod3} performs better than Model~\ref{mod2}. Therefore, we restrict the following tests to  Model~\ref{mod3}. We compare the impact of employing different objective functions from Section~\ref{sec:obj} on the economic efficiency and the customer satisfaction of the final routing solution. The respective objective functions are used in Model~\ref{mod3} in the place of \eqref{costobj}. Moreover, for the objective function $f_{rcr}$ we add variables $p_i$, $i\in R$ to Model~\ref{mod3} and replace constraints \eqref{passpickedup1} by constraints \eqref{passpick?}. In case of the objective functions $f_r$, $f_{cr}$ and $f_{rcr}$, we add variables $d_i \geq 0$, $i\in R$ and constraints \eqref{regret} to Model~\ref{mod3}.  In case of the objective functions involving the users' \emph{maximum regret}, i.e.\ $f_{r_{\max}}$ and $f_{cr_{\max}}$, we additionally add the variable $d_{\max} \geq 0$ and constraints \eqref{maxregret} to Model~\ref{mod3}.

The weights in the objective functions involving more than one criterion are chosen from a user perspective and based on preliminary numerical tests. In the first part of the computational experiments we consider the single-objective functions $f_c$, $f_r$ and $f_{r_{\max}}$. The weights in $f_{cr}$ and $f_{cr_{\max}}$  are then chosen so that the values of total and maximum regret, respectively,  in $f_{cr}$ and $f_{cr_{\max}}$ deviate not more than 15\% from the optimal objective values for $f_r$ and $f_{r_{\max}}$ (in instances solved to optimality). After some preliminary testing the weights are set to $\alpha=1$ and $\beta=\frac{n}{5}$, which yields
\begin{align*}
	f_{cr} & = \sum_{a \in A}c_a \, x_a +  \sum_{i\in R}d_i \quad \text{and} \\
	f_{cr_{\max}} & = \sum_{a \in A}c_a \, x_a + \frac{n}{5}\, d_{\max}.
\end{align*}
Note that choosing the weighting parameter $\beta$ as a function of $n$ ensures that not only the first term in $f_{cr_{\max}}$ grows with the number of users.

By using the objective function $f_{rcr}$ we can measure the impact of denying certain requests on the served users' regret and the total routing costs. The weighting parameter $\gamma$ in $f_{rcr}$ defines the trade-off between the general requirement of answering as many requests as possible on one hand, and the goal of cost and time efficient transport solutions on the other hand. After testing several weights, it turns out that $\gamma = 20$ is a reasonable choice, yielding 
\begin{equation*}
	f_{rcr} = \sum_{a \in A} c_a \, x_a  + \sum_{i \in R} d_i + 20\left(n - \sum_{i \in R} p_i \right). 
\end{equation*}
In our numerical experiments, on average at most 10\% of the requests are rejected when using these weights. Note that several authors using weighted sum objectives base their choice of weights on \cite{Jorgensen2007} (see e.g.\ \cite{Mauri2009, Kirchler2013}). However, the total routing costs, the total/maximum regret and the number of unanswered requests depend strongly on the test instances and may vary considerably. Since in this section we create a new class of test instances, we refrain from using these predetermined weights. 

\begin{table}
	\centering
	\begin{threeparttable}
		\footnotesize
		\begin{tabular}{cccccccc}
			\toprule 
			& \multicolumn{3}{c} {$f_c$} & \multicolumn{4}{c}{$f_r$}\\
			\cmidrule(lr){2-4}\cmidrule(lr){5-8}\addlinespace[0.1em]
			{$n$} & {$f_c$}	& {$f_r$}	& {CPU} & {$f_c$} & {$f_r$} & {Gap} & {CPU} \\
			\midrule 
			\textbf{Q3} &  &  &  &  &  &  & \\ \addlinespace[0.2em]
			20 & 113.3 & 136.8 & 0.06 & 142.5 & 12.2 & 0 & 0.05 \\
			30 & 170.7 & 250.3 & 0.06 & 209.1 & 23.5 & 0 & 0.07 \\
			40 & 205.8 & 346.6 & 0.15 & 259.3 & 46.4 & 0 & 0.55 \\
			60 & 309.6 & 631.0 & 1.09 & 385.7 & 128.8 & 0 & 44.33 \\
			80 & 383.5 & 798.0 & 7.00 & 480.2\tnote{4} & 193.6\tnote{4} & 13\tnote{4} & 4721\tnote{4} \\
			100 & 445.5 & 1070.4 & 218 & 603.7\tnote{3} & 193.0\tnote{3}& 31\tnote{3} & 6060\tnote{3} \\
			\midrule
			\textbf{Q6} &  &  &  &  &  &  & \\ \addlinespace[0.2em]
			20 & 117.6 & 144.2 & 0.02 & 135.7 & 14.3 & 0 & 0.02 \\
			30 & 193.5 & 202.3 & 0.04 & 223.7 & 39.8 & 0 & 0.05 \\
			40 & 246.8 & 273.0 & 0.06 & 294.8 & 48.16 & 0 & 0.06 \\
			60 & 343.9 & 467.3 & 0.26 & 393.3 & 99.7 & 0 & 2.13 \\
			80 & 448.9 & 640.4 & 0.73 & 520.3 & 148.8 & 0 & 106 \\
			100 & 543.4 & 791.6 & 6.07 & 627.1 & 200.7 & 5 & 3236 \\
			\bottomrule \addlinespace[0.2em]
		\end{tabular}
		\begin{tablenotes}
			\small
			\item[3] Average has been built only over three instances, as two instances were not solved within 7200s.
			\item[4] Average has been built only over four instances, as one instance was not solved within 7200s.
		\end{tablenotes}
		\caption{Average values on the Q3 and Q6 test instances solved with the objective functions $f_c$ and $f_r$.}
		\label{tab:cr}
	\end{threeparttable}
\end{table}

\begin{table}
	\centering
	\begin{threeparttable}
		\footnotesize
		\begin{tabular}{cccc@{\hspace{1.2pt}}c@{\hspace{1.2pt}}ccccc@{\hspace{1.2pt}}c@{\hspace{1.2pt}}c}
			\toprule 
			& \multicolumn{5}{c}{$f_{cr}$} & \multicolumn{6}{c}{$f_{rcr}$} \\
			\cmidrule(lr){2-6}\cmidrule{7-12}\addlinespace[0.1em]
			{$n$} & {Obj.v.} & {$f_c$} & {$f_r$} & {Gap} & {CPU}  & {Obj.v.}& {$f_c$} & {$f_r$} & {a.r.} & {Gap} & {CPU} \\
			\midrule
			\textbf{Q3}& & & & & & & & && & \\ \addlinespace[0.2em]
			20 & 141.5 & 128.3 & 13.2 & 0 & 0.04 & 140.4 & 121.2 & 7.2 & 97 & 0 & 0.04 \\
			30 & 219.0 & 192.5 & 26.4 & 0 & 0.08 & 215.9 & 186.5 & 17.3 & 98 & 0 & 0.08 \\
			40 & 284.0 & 234.3 & 49.7 & 0 & 0.47 & 275.3 & 220.1 & 19.2 & 96 & 0 & 0.47 \\
			60 & 495.6 & 361.6 & 134.0 & 0 & 37.04 & 459.3 & 326.7 & 56.6 & 93 & 0 & 6.77 \\
			80 & 663.3\tnote{4} & 454.3\tnote{4} & 208.9\tnote{4} & 3\tnote{4} & 3680\tnote{4}& 598.2 & 412.8 & 61.4 & 92 & 1.2 & 1705 \\
			100 & 762.7 & 533.8 & 228.9 & 12 & 6416 & 718.6 & 500.3 & 86.4 & 93 & 8.2 & 5796 \\
			\midrule
			\textbf{Q6}& & & & & & && & & & \\ \addlinespace[0.2em]
			20 & 139.8 & 124.5 & 15.4 & 0 & 0.03 & 138.1 & 122.5 & 11.6 & 99 & 0 & 0.03 \\
			30 & 253.6 & 209.2 & 44.4 & 0 & 0.07 & 241.5 & 189.8 & 11.7 & 93 & 0 & 0.05 \\
			40 & 316.3 & 267.4 & 48.9 & 0 & 0.08 & 304.1 & 249.0 & 19.1 & 96 & 0 & 0.07 \\
			60 & 476.9 & 373.7 & 103.2 & 0 & 1.52 & 438.1 & 340.7 & 33.4 & 95 & 0 & 0.64 \\
			80 & 655.4 & 503.3 & 152.1 & 0 & 20.24 & 589.5 & 450.1 & 27.4 & 93 & 0 & 1.75 \\
			100 & 804.9 & 595.0 & 209.9 & 0 & 853 & 719.5 & 530.7 & 44.8 & 93 & 0 & 5.28 \\
			\bottomrule\addlinespace[0.2em]
		\end{tabular}
		\begin{tablenotes}
			\small
			\item[4] Average has been built only over four instances, as one instance was not solved within 7200s.
		\end{tablenotes}
		\caption{Average values on the Q3 and Q6  test instances solved with the objective function $f_{cr}$ and $f_{rcr}$.}
		\label{tab:crrcr}
	\end{threeparttable}
\end{table}

\begin{table}
	\centering
	\begin{threeparttable}
		\footnotesize
		\begin{tabular}{cccccccccccc}
			\toprule 
			& \multicolumn{5}{c}{$f_{r_{\max}}$}& \multicolumn{6}{c}{$f_{cr_{\max}}$} \\
			\cmidrule(lr){2-6}\cmidrule(lr){7-12}\addlinespace[0.1em]
			{$n$} &{$f_c$}   &    {$f_r$}  &{$f_{r_{\max}}$}&Gap&CPU&{Obj.v.}&{$f_c$} &{$f_r$} &{$f_{r_{\max}}$}&Gap&CPU\\
			\midrule
			\textbf{Q3}& & & & & & & & & & & \\ \addlinespace[0.2em]
			20 & 142.0 & 48.6 & 8.9 & 0 & 0.03 & 159.0 & 119.9 & 91.9 & 9.8 & 0 & 0.04 \\
			30 & 211.7 & 64.6 & 7.9 & 0 & 0.05 & 239.1 & 191.8 & 114.3 & 7.9 & 0 & 0.14 \\
			40 & 255.9 & 132.5 & 9.5 & 0 & 0.16 & 305.8 & 228.4 & 218.4 & 9.7 & 0 & 0.65 \\
			60 & 388.8 & 301.4 & 13.0 & 0 & 0.84 & 484.6 & 325.6 & 442.9 & 13.3 & 0 & 3.11 \\
			80 & 495.0 & 427.5 & 12.6 & 0 & 1443 & 634.2 & 431.3 & 557.4 & 12.7 & 1 & 1749 \\
			100 & 599.9 & 499.2 & 11.0 & 1 & 1496 & 738.0 & 506.5 & 709.7 & 11.6 & 2 & 5884 \\
			\midrule
			\textbf{Q6}& & & & & & & & & & & \\ \addlinespace[0.2em]
			20 & 137.4 & 42.2 & 5.9 & 0 & 0.02 & 149.8 & 125.9 & 50.6 & 6.0 & 0 & 0.04 \\
			30 & 223.8 & 118.6 & 11.2 & 0 & 0.11 & 274.3 & 206.8 & 149.9 & 11.2 & 0 & 0.05 \\
			40 & 293.5 & 114.0 & 12.3 & 0 & 0.06 & 355.9 & 255.0 & 232.5 & 12.6 & 0 & 0.11 \\
			60 & 406.8 & 246.7 & 11.5 & 0 & 0.28 & 502.4 & 363.8 & 332.1 & 11.6 & 0 & 0.84 \\
			80 & 531.2 & 365.9 & 11.0 & 0 & 1.52 & 663.0 & 485.7 & 431.9 & 11.1 & 0 & 7.15 \\
			100 & 646.3 & 478.8 & 12.0 & 0 & 15.93 & 815.3 & 573.2 & 588.3 & 12.1 & 0 & 71.12 \\
			\bottomrule \addlinespace[0.2em]
		\end{tabular}
		\caption{Average values on the Q3 and Q6 test instances solved with the objective functions $f_{r_{\max}}$ and $f_{cr_{\max}}$.}
		\label{tab:rmaxcrmax}
	\end{threeparttable}
\end{table}

In Tables~\ref{tab:cr}, \ref{tab:crrcr} and \ref{tab:rmaxcrmax} the results are summarized. All figures reported are average values over five instances $Q3.n.m$, $m\in \{1,\ldots,5\}$ and $Q6.n.m$, $m\in \{1,\ldots,5\}$, respectively. The tables contain information about the total routing  costs, the total regret (and where applicable the maximum regret and the percentage of answered requests) for each number of user requests $n\in\{20,30,40,60,80,100\}$ of the instance sets $Q3$ and $Q6$. Furthermore, the relative gap and the CPU time returned by CPLEX are reported. If no relative gap is shown, all instances have been solved to optimality.

The computational times show that Model~II combined with a cost-objective function can be used to solve small to medium-sized artificial instances very efficiently. This is in acccordance with the results indicated by the benchmark instances. The computational time increases when other objective functions are used, but most of the instances can still be solved within a few seconds. However, for the larger $Q3$-instances ($n\in \{80,100\}$) the CPU time increases drastically. 
In general, the $Q6$-instances are easier to solve than the $Q3$-instances, reflecting the smaller number of possible user allocations in the vehicle when the number of requested seats per user is not limited to one. 
While the scenario that users may request any number of seats between one and six reflects the conditions under which ridepooling services operate in reality, a uniform distribution of $q_i$ in the set $\{1, \ldots, 6\}$ is probably not realistic. In future work, this should be tested at real world data. 

\begin{table}
	\centering
	\begin{threeparttable}
		\footnotesize
		\begin{tabular}{ccccccccccccc}
			\toprule 
			& \multicolumn{2}{c} {$f_r$ vs. $f_c$ } & \multicolumn{2}{c}{$f_{cr}$ vs. $f_c$}& \multicolumn{2}{c}{$f_{cr}$ vs. $f_r$} & \multicolumn{2}{c}{$f_{r_{\max}}$ vs. $f_c$} &  \multicolumn{2}{c}{$f_{cr_{\max}}$ vs. $f_c$} & \multicolumn{2}{c}{$f_{cr_{\max}}$ vs. $f_{r_{\max}}$}\\
			\cmidrule(lr){2-3}\cmidrule(lr){4-5}\cmidrule(lr){6-7}\cmidrule(lr){8-9}\cmidrule(lr){10-11}\cmidrule(lr){12-13}\addlinespace[0.1em]
			$n$ & $\Delta f_c$	& $\Delta f_r$	& $\Delta f_c$  & $\Delta \text{CPU}$ &  $\Delta f_r$ &  $\Delta \text{CPU}$ &  $\Delta f_c$ & $\Delta f_r$ & $\Delta f_c$  & $\Delta$ CPU &  $\Delta f_{r_{\max}}$ &  $\Delta$ CPU \\
			\midrule
			\textbf{Q3} & & & & & & & & & & & & \\ \addlinespace[0.2em]
			20 & 26 & -91 & 13 & -36 & 8 & -22 & 25 & -64 & 6 & -21 & 10 & 47 \\
			30 & 22 & -91 & 13 & 32 & 12 & 24 & 24 & -74 & 12 & 119 & 0 & 152 \\
			40 & 26 & -87 & 14 & 216 & 7 & -15 & 24 & -62 & 11 & 336 & 2 & 309 \\
			60 & 25 & -80 & 17 & 3286 & 4 & -16 & 26 & -52 & 5 & 184 & 2 & 269 \\
			80 & 25 & -76 & 18 & 52539 & 8 & -22 & 29 & -46 & 12 & 24927 & 1 & 21 \\
			100 & 35 & -82 & 20 & 2840 & 19 & 6 & 35 & -53 & 14 & 2597 & 6 & 293 \\
			Avg. & 27 & -84 & 16 & 9813 & 10 & -7 & 27 & -59 & 10 & 4690 & 4 & 182 \\
			\midrule
			\textbf{Q6} & & & & & & & & & & & &\\ \addlinespace[0.2em]
			20 & 15 & -90 & 6 & 75 & 7 & 75 & 17 & -71 & 7 & 150 & 1 & 82 \\
			30 & 16 & -80 & 8 & 57 & 12 & 43 & 16 & -41 & 7 & 29 & 0 & -51 \\
			40 & 19 & -82 & 8 & 36 & 2 & 23 & 19 & -58 & 3 & 96 & 3 & 77 \\
			60 & 14 & -79 & 9 & 489 & 4 & -29 & 18 & -47 & 6 & 224 & 1 & 203 \\
			80 & 16 & -77 & 12 & 2664 & 2 & -81 & 18 & -43 & 8 & 877 & 0 & 372 \\
			100 & 15 & -75 & 9 & 13948 & 5 & -74 & 19 & -40 & 5 & 1071 & 1 & 346 \\
			Avg. & 16 & -80 & 9 & 2878 & 5 & -7 & 18 & -50 & 6 & 408 & 1 & 172 \\
			\bottomrule \addlinespace[0.2em]
		\end{tabular}
		\caption{Comparison of the change in total routing costs, total regret and CPU time (all in percent) on the Q3 and Q6 test instances solved with the objective functions $f_c$, $f_r$, $f_{cr}$, $f_{r_{\max}}$ and $f_{cr_{\max}}$.}
		\label{tab:compcrcrrcr}
	\end{threeparttable}
\end{table}

\begin{table}
	\centering
	\begin{threeparttable}
		\footnotesize
		\begin{tabular}{cccccccc}
			\toprule 
			& \multicolumn{3}{c} {$f_{cr_{\max}}$ vs. $f_{cr}$ } & \multicolumn{4}{c}{$f_{rcr}$ vs. $f_{cr}$}\\
			\cmidrule(lr){2-4}\cmidrule(lr){5-8}\addlinespace[0.1em]
			$n$ & $\Delta f_c$	& $\Delta f_r$ & $\Delta \text{CPU}^\star$ &  $\Delta f_c$ & $\Delta f_r$ & $\Delta \text{CPU}^\star$& $\Delta $ a.r. \\
			\midrule
			\textbf{Q3} & & & & & & & \\ \addlinespace[0.2em]
			20 & -7 & 594 & 22 & -6 & -46 & 6 & -3 \\
			30 & 0 & 333 & 66 & -3 & -34 & 0 & -2 \\
			40 & -3 & 340 & 38 & -6 & -61 & 1 & -4 \\
			60 & -10 & 231 & -92 & -10 & -58 & -82 & -7 \\
			80 & -5 & 167 & -52 & -9 & -71 & -54 & -8 \\
			100 & -5 & 210 & -8 & -6 & -62 & -10 & -7 \\
			Avg  & -5 & 312 & -4 & -7 & -55 & -23 & -5 \\
			\midrule
			\textbf{Q6}  & & & & & & &  \\ \addlinespace[0.2em]
			20 & 1 & 230 & 43 & -2 & -24 & 21 & -1 \\
			30 & -1 & 237 & -18 & -9 & -74 & -21 & -7 \\
			40 & -5 & 375 & 45 & -7 & -61 & -5 & -4 \\
			60 & -3 & 222 & -45 & -9 & -68 & -58 & -5\\
			80 & -3 & 184 & -65 & -11 & -82 & -91 & -7 \\
			100 & -4 & 180 & -92 & -11 & -79 & -99 & -7 \\
			Avg  & -2 & 238 & -22 & -8 & -65 & -42 & -5 \\
			\bottomrule \addlinespace[0.2em]
		\end{tabular}
		\caption{Comparison of the change in total routing costs, total regret, CPU time and number of answered requests (all in percent) on the Q3 and Q6 test instances solved with the objective functions $f_{cr}$, $f_{cr_{\max}}$ and $f_{rcr}$.}
		\label{tab:crcrmaxrcr}
	\end{threeparttable}
\end{table}


A comparison of the effects of different objective functions in Model~\ref{mod3}  can be found in Tables~\ref{tab:compcrcrrcr} and \ref{tab:crcrmaxrcr}. In the second and third column of Table~\ref{tab:compcrcrrcr} we illustrate the change (in percent) in total routing costs and total regret when replacing $f_c$ by $f_r$. While the regret decreases by an average of 84\% and 80\% (Q3- and Q6-instances, respectively), the routing costs only increase by 27\% and 16\% on average. This shows that by including the criterion of the users' regret in the objective function we can spare users a large amount of unnecessary ride or waiting time. This comes at the cost of higher routing costs. However, even from a service provider's perspective it might be reasonable to accept a rather small loss in order to improve user convenience and to remain competitive. In column six of the same table we demonstrate that the weights chosen in the cost-regret objective $f_{cr}$ indeed reflect a user perspective: 
On average (over the instances solved to optimality within the time limit) there is an increase of at most 15\%  in total regret compared to solely optimizing w.r.t.\ $f_r$. Moreover, there is an average increase in routing costs of 16\% (Q3-instances) and 9\% (Q6-instances) compared to the costs when using routing costs as the only optimization criterion, i.e.\ when using $f_c$ as the objective function. Comparing $f_{cr}$ to $f_c$,  we observe an increase in CPU time for all instances but $Q=3$ and $n=20$, but CPU times for $Q=3$ and $n\in\{20,30,40,60\}$ and all Q6-instances are still reasonable as shown in Table~\ref{tab:crrcr}.  Comparing $f_{cr}$ to $f_r$ we observe an average decrease in CPU time of 7\% for the Q3- and Q6-instances. Similar results are obtained for the objective functions $f_c$, $f_{r_{\max}}$ and $f_{cr_{\max}}$, although the average decrease in regret when comparing $f_{r_{\max}}$ to $f_c$ is only 59\% (Q3-instances) and 50\% (Q6-instances). The meaningful choice of the weighting parameter in $f_{cr_{\max}}$ is demonstrated by the second last column in Table~\ref{tab:compcrcrrcr}, as the increase in maximum regret is below 15\% for all instances. 

Since both of the weighted sum objectives $f_{cr}$ and $f_{cr_{\max}}$ improve user convenience and increase routing costs compared to $f_c$, we evaluate which of the objective functions is more effective in this respect. Columns two to four in Table~\ref{tab:crcrmaxrcr} illustrate that when using $f_{cr_{\max}}$ instead of $f_{cr}$ the average computational time decreases between 8\% to 92\% (Q3-instances) and between 45\% and 92\% (Q6-instances) in the larger instances, i.e.\ $n\in\{60,80,100\}$, which are harder to solve in general. Despite its computational superiority and the fact that the objective function $f_{cr_{\max}}$ generally improves the user satisfaction of optimal tours, it has some shortcomings. Indeed, if there is one user with a high regret $d_i$, for instance, because he or she is the last user in a tour to be dropped-off, the time loss of all other users $j$ with $d_j < d_i$ has no impact on the objective function. Therefore, the other users may not be driven to their drop-off points as quickly as possible. This is reflected by the increase in regret using the objective $f_{cr_{\max}}$ as compared to $f_{cr}$, shown in Table~\ref{tab:crcrmaxrcr}: 312\% (Q3-instances) and 238\% (Q6-instances). The routing costs remain roughly the same; there is an average decrease of 5\% (Q3-instances) and 2\% (Q6-instances). Due to these shortcomings, we do not consider a tri-criterion weighted sum objective function $f_{rcr_{\max}}$, but restrict our attention to $f_{rcr}$ in the remaining discussion.

\begin{figure}
	\centering
	\includegraphics[width=0.7\textwidth, clip,bb= 50 43 392 338]{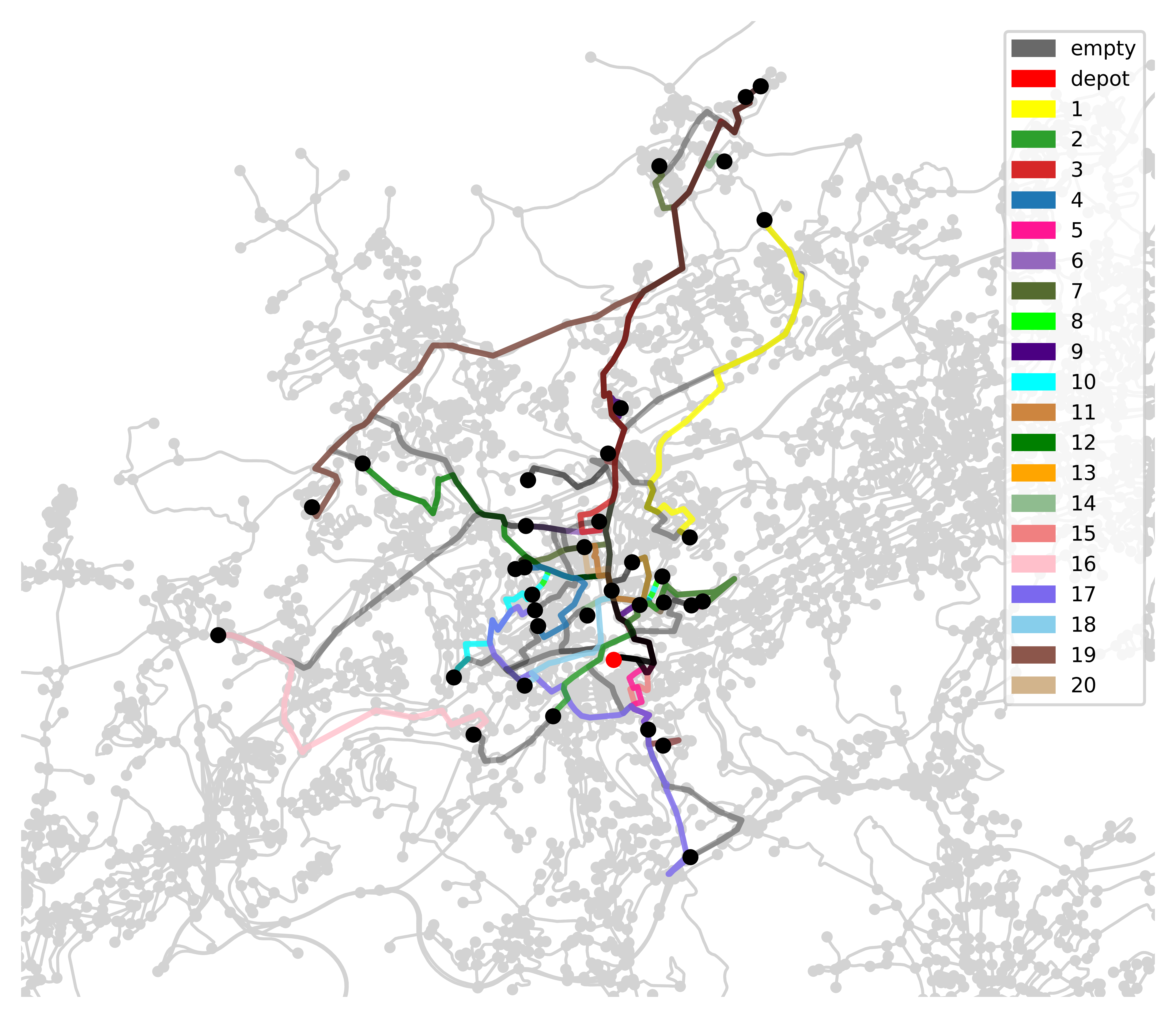}
	\caption{Vehicle routes of instance $Q3n20.5$ solved using the objective function $f_{cr}$.}
	\label{figce}
\end{figure}

\begin{figure}
	\centering
	\includegraphics[width=0.7\textwidth, clip,bb= 50 43 392 338]{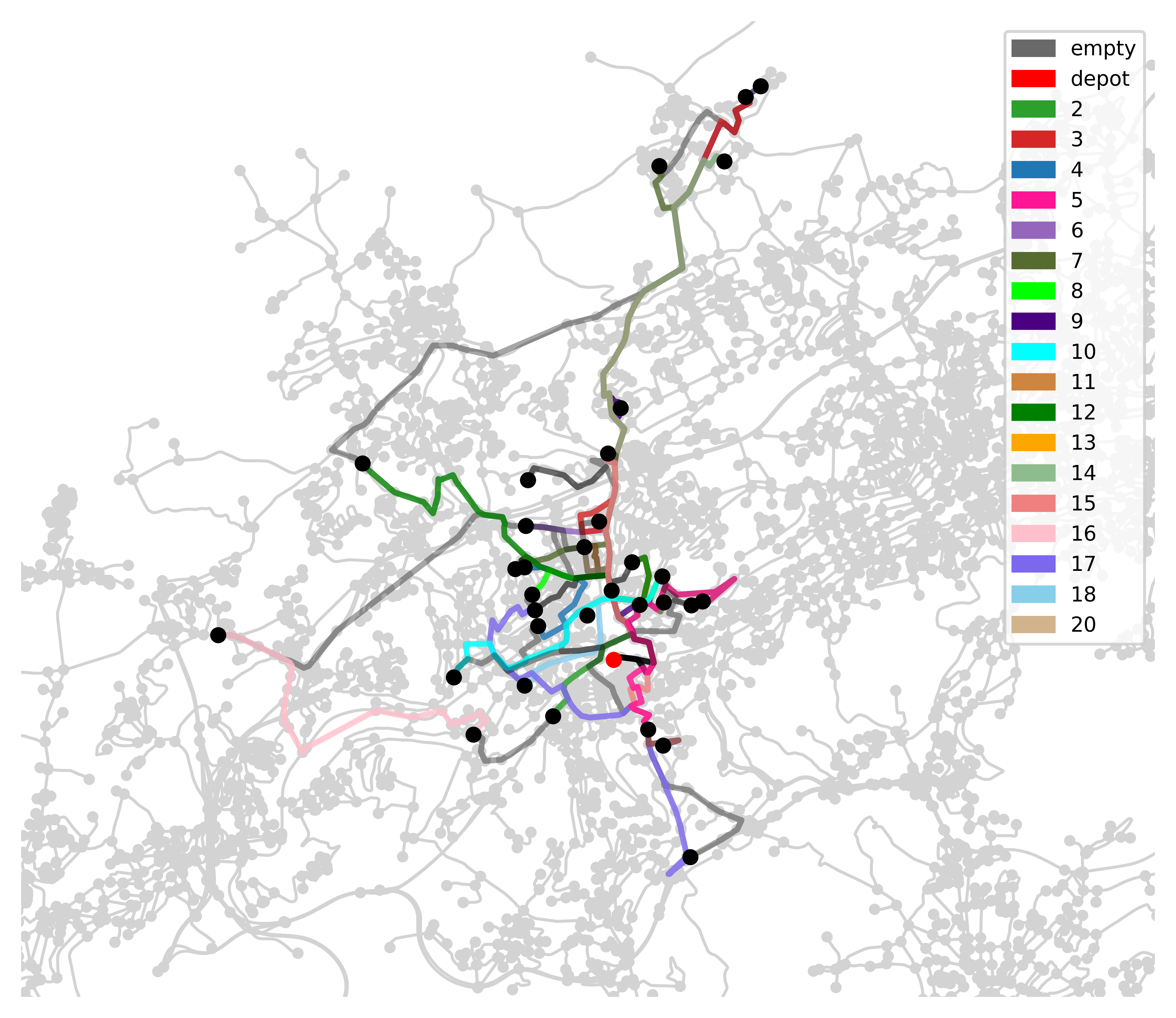}
	\caption{Vehicle routes of instance $Q3n20.5$ solved using the objective function $f_{rcr}$.}
	\label{figrce}
\end{figure}

\begin{table}
	\begin{threeparttable}
		\footnotesize
		\begin{tabularx}{\textwidth}{p{1cm}p{1cm}XXXXXXXXXXXX}
			\toprule\addlinespace[0.2em]
			\multicolumn{14}{c}{$f_{cr}$} \\
			\midrule\addlinespace[0.5em]
			\textbf{Tour 1\quad}
			& Location & $1^+$ & $1^-$ & $15^+$ & $15^-$ & $5^+$ & $5^-$ & $8^+$ & $10^+$ & $8^-$ & $10^-$ & & \\
			&Time[m] & 24.4       & 43.2      & 60.0     & 76.9   & 120.0  & 136.6 & 156.5  & 160.0 & 168.0 & 175.6 & & \\[2ex]
			\textbf{Tour 2}
			&Location & $9^+$ & $9^-$ & $6^+$ & $6^-$ & $7^+$ & $7^-$ & & & & & &\\
			&Time[m]  & 60.0    & 71.7   & 95.0    & 120.6 & 150.0  & 173.2& &  & & & &\\[2ex]
			\textbf{Tour 3}
			&Location  & $12^+$ & $12^-$ & $14^+$ & $14^-$& $17^+$ & $17^-$ & & & & & &\\
			&Time[m]   & 60.0     & 79.4      & 105.0     & 128.0  & 180.0    & 200.6    & & & & & &  \\[2ex]
			\textbf{Tour 4}
			&Location   & $13^+$ & $13^-$ & $3^+$ & $3^-$ & $19^+$ & $19^-$ &  &  &  &  & &\\
			&Time[m]  & 30.0      & 31.0      & 85.0     &108.9 & 140.0    & 170.7    &  & &  &  & &\\[2ex]
			\textbf{Tour 5}
			&Location   & $20^+$ & $20^-$ & $4^+$ & $4^-$ & $18^+$ & $18^-$ & $11^+$ & $11^-$ & $16^+$& $16^-$ & $2^+$ & $2^-$ \\
			&Time[m]  & 50.0       & 56.2     & 65.0    & 72.0    & 85.0       & 92.5     & 95.0     & 102.4  & 121.3 & 138.1 & 155.0 & 169.1\\ 
			\bottomrule
		\end{tabularx}
		\caption{Vehicle routes (without depot) of instance $Q3n20.5$ solved using the objective function $f_{cr}$.} 
		\label{tab:figce}
	\end{threeparttable}
\end{table}

\begin{table}
	\begin{threeparttable}
		\footnotesize
		\begin{tabularx}{\textwidth}{p{1cm}p{1cm}XXXXXXXXXXXX}
			\toprule\addlinespace[0.2em]
			\multicolumn{14}{c}{$f_{rcr}$} \\
			\midrule\addlinespace[0.5em]
			\textbf{Tour 1\quad}
			& Location & $15^+$ & $15^-$ & $5^+$ & $5^-$ & $10^+$ & $10^-$ &  &  & & & &\\
			&Time[m]  & 60.0       & 76.9      & 120.0  & 135.6 & 160.0    & 172.5   &  &  & & & &\\[2ex]
			\textbf{Tour 2}
			&Location  & $13^+$ & $13^-$ & $9^+$ & $9^-$ & $6^+$ & $6^-$ & $7^+$ & $7^-$ & & & &\\
			&Time[m]  & 30.0      & 31.0      & 60.0    & 71.7   & 95.0    & 120.6 & 150.0  & 173.2& & & &\\[2ex]
			\textbf{Tour 3}
			&Location  & $12^+$ & $12^-$ & $14^+$ & $14^-$& $17^+$ & $17^-$ & & & & & &\\
			&Time[m]  & 60.0     & 79.4      & 105.0     & 128.0  & 180.0    & 200.6    & & & &  & & \\[2ex]
			\textbf{Tour 4}
			&Location   & $3^+$ & $3^-$ & $8^+$ & $8^-$ &  &  &  &  &  & & & \\
			&Time[m]  & 85. 0    & 108.9 & 155.0 & 163.0 &  &  &  &  &  & & & \\[2ex]
			\textbf{Tour 5}
			&Location   & $20^+$ & $20^-$ & $4^+$ & $4^-$ & $18^+$ & $18^-$ & $11^+$ & $11^-$ & $16^+$& $16^-$ & $2^+$ & $2^-$ \\
			&Time[m]  & 50.0       & 56.2     & 65.0    & 72.0    & 85.0       & 92.5     & 95.0     & 102.4  & 121.3 & 138.1 & 155.0 & 169.1\\ 
			\bottomrule
		\end{tabularx}
		\caption{Vehicle routes (without depot) of instance $Q3n20.5$ solved using the objective function $f_{rcr}$.} 
		\label{tab:figrce}
	\end{threeparttable}
\end{table}


The last four columns of Table~\ref{tab:crcrmaxrcr}, $f_{rcr}$ vs. $f_{cr}$, illustrate the decrease in the overall routing costs and regret that is obtained when rejecting some of the requests. This is illustrated at the instance $Q3n20.5$ in Figure~\ref{figce}, showing an optimal tour w.r.t.\ $f_{cr}$, and Figure~\ref{figrce} where an optimal tour w.r.t.\ $f_{rcr}$ is shown. In Figure~\ref{figce} the driver has to make a detour to transport requests $1$ and $9$, while it becomes obvious from Figure~\ref{figrce} that the service provider benefits from a large decrease in routing costs. Table~\ref{tab:figce} contains the corresponding vehicle routes including pick-up and drop-off times in minutes after the start of service. In the vehicle tours computed using $f_{rcr}$, users $8$ and $10$ who have been transported with an associated regret of 5 and 3 minutes, respectively, are transported without any time loss.

For the Q3- and Q6-instances, on average 5\% less of the requests are answered in comparison to the results obtained with the objective function $f_{cr}$. In turn, we observe an average decrease of the total routing costs and the total regret for all instance sizes. While the decrease in routing costs ranges from 2\% to 11\%, there are huge savings in regret: There is an average decrease of 55\% and 65\% (Q3- and Q6-instances, respectively.)  Furthermore, computational time decreases by 23\% (Q3-intances) and 42\% (Q6-instances). 

While it comes at cost of an increase in computational time and routing costs, the tests show that it is worthwhile to use more than one optimization criterion. Objective functions $f_{cr}$ and $f_{rcr}$ provide satisfactory results w.r.t.\ an increase in the users' total regret as compared to the total regret when using it as the only optimization criterion. Moreover, we have shown that the service quality can be improved for some users and costs can be reduced by rejecting unfavourable requests. Average computation times using objective function $f_{rcr}$ for $Q=3$, $n\in\{20,30,40,60\}$ and $Q=6$ are less than 7 seconds. Hence, our model can be used to solve small to medium-sized instances. It is applicable without the need for extensive coding, as it only consists of an MIP solver and the generation of the event-based graph.

\section{Conclusions}\label{sec:concl}
\noindent
In this paper we suggest a new perspective on modeling ridepooling problems. By using an event-based graph representation rather than a geographical model, we show that capacity, pairing and precedence constraints can be handled implicitly. While the resulting MILP formulations are pseudo-compact and generally have more variables as compared to classical compact models, extensive numerical experiments show that the implicit constraint formulation leads to considerably improved computational times. Indeed, both for benchmark instances from the literature as well as for artificial instances in the city of Wuppertal, problems with up to $80$ requests can be solved in less than 7 seconds. 

Moreover, we analyse the effects of including additional optimization criteria in the model. In addition to the classical cost objective function, we consider the (total or maximum) regret as well as the number of rejected requests as measures for user satisfaction. By combining these criteria into a weighted sum objective function we demonstrate that user satisfaction can be largely improved at only relatively small additional expenses, i.e., overall routing costs.

While most real world instances of the DARPs are much larger, the proposed approach could still prove useful as a subroutine also for realistically-sized instances. Particularly, the extension to the online version of the DARP in which requests arrive over time could be a promising application for our models, see \cite{Gaul2021A}. Indeed, often only very few new requests arrive simultaneously and re-routing of already scheduled users is only acceptable if it decreases their arrival time. Consequently, the number of simultaneous users in such a rolling-horizon version of a dial-a-ride problem is relatively small and could potentially be solved exactly using one of our models.
Furthermore, we will investigate the integration of ridepooling services into the local public transport system in Wuppertal.

\section*{Acknowledgements}
\noindent This work was partially supported by the state of North Rhine-Westphalia (Germany) within the project ``bergisch.smart.mobility''. Furthermore, we thank three anonymous reviewers for their valuable comments.

\newpage

\bibliography{ridehailing}

\end{document}